\documentclass[12pt,oneside]{amsart}
       \usepackage{amssymb}
       \usepackage{epsfig}
       \usepackage{amsfonts}
       \usepackage{eucal}
       \usepackage{amsmath,amsthm,amsfonts,latexsym,amscd,epsfig}
       \usepackage[notref,notcite]{}
       \usepackage{pictex}
       \theoremstyle{plain}
       \newtheorem{theorem}{Theorem}[section]
       \newtheorem{lemma}[theorem]{Lemma}
       \newtheorem{corollary}[theorem]{Corollary}
       
       \newtheorem{proposition}[theorem]{Proposition}
       
       \newtheorem{definition-theorem}[theorem]{Definition}

      \newtheorem{definition}{Definition}[section]
      \newtheorem{remark}{Remark}[section]

       \newenvironment{defn}{\begin{definition}\begin{rm}}
                {\end{rm}\end{definition}}

       \newcommand{\R}{{\mathbb R}}
       \newcommand{\C}{{\mathbb C}}
       
       \newcommand{\Z}{{\mathbb Z}}

       \newcommand{\grtimes}{\hat{\otimes}}
       
       \newcommand{\del} {\partial}

       \newcommand{\alg}{{\mathbf a}}

       \def\co{\colon\thinspace}
       \makeatletter
       \def\@setcopyright{}
       \def\serieslogo@{}
       \makeatother

\begin{document}

    \author{Dmitry Gerenrot}
    \address{ Dmitry Gerenrot, Georgia Institute of Technology,
    School of Mathematics, 686 Cherry Street, Atlanta,
    GA, 30332-0160, U.S.A.}
    \email{gerenrot@math.gatech.edu}


\thanks{ The author expresses his most sincere gratitude to his thesis
advisor,  professor Nigel Higson for invaluable guidance, patience
and summer support.}

    \title[Residue Formulation]{Residue Formulation of
        Chern Character on Smooth Manifolds}

\begin{abstract}
The Chern character of a complex vector bundle is
most conveniently defined as the exponential
of a curvature of a connection.  It is well known that its cohomology
class does not depend on the  particular
connection chosen. It has been shown by Quillen \cite{Q}
that a connection may be perturbed by an endomorphism
of the vector bundle, such as a symbol of some elliptic
differential operator. This point of view, as we intend to show,
allows one to relate Chern character to a non-commutative
sibling formulated by Connes and Moscovici \cite{CM}.
The general setup for our problem is purely geometric.
Let $\sigma$ be the symbol of a Dirac-type operator acting
on sections of a ${\Z}_2$-graded vector bundle $E$.
Let $\nabla$ be a connection on $E$, pulled back to
$T^*M$. Suppose also that $\nabla$ respects the
${\Z}_2$-grading. The object $\nabla+\sigma$ is a superconnection
on $T^*M$ in the sense of Quillen.  We obtain a formula for the
$H_*(M)$-valued Poincare dual of Quillen's Chern character
$ch(D)=\operatorname{tr_s}e^{(\nabla+\sigma)^2}$  in terms of residues of
$\Gamma(z)\operatorname{tr_s}(\nabla+\sigma)^{-2z}$.  We also
compute two examples.

\end{abstract}


    \date{\today}

    \maketitle

\section{Introduction.}\label{section:intro}

The historical background for noncommutative index theory has two
basic parts. The first one dates back to the nineteen fifties, when Israel
Gelfand pointed out to Sir Michael Atiyah that the index of a
Fredholm operator was stable under small perturbations. The
ultimate consequence of this remark is quite famous: the Atiyah-Singer
Index Theorem \cite{AS0, AS123, AS4}. The second part is development of noncommutative geometry by
Alain Connes \cite{Con,GVF}.

In particular, two noncommutative versions of the Chern  character were
developed. There is one due to Connes \cite{NCDG} and a more
recent one due to  Connes and Moscovici \cite{CM}. We are interested in
this more recent version, which is called the {\it residue cocycle}. Its
individual terms are certain residues which are geometrically
interesting, as pointed out by Higson. However, they are not very
well understood.

In the present paper, we shall  prove a formula which
resembles the formula of Connes and Moscovici in \cite{CM}, albeit
in a classical setting.
Suppose $M$ is a compact smooth manifold with no boundary.
Let $E\rightarrow M$ be a ${\Z}_2$-graded complex vector bundle.
Let $D$ be an elliptic, odd, first-order, skew-adjoint
differential operator on $E$. Finally, let $\pi\co T^*M\to M$ be
the standard projection map of the cotangent bundle.  With this setup,
we will establish a formula for the Chern character of the symbol of $D$
which resembles the Connes-Moscovici Chern character in noncommutative
geometry. It is comprised of finitely many residues of zeta functions
constructed from the symbol of $D$ and a connection on $E$.

We shall use Quillen's formalism in which the symbol $L$ of $D$,
together with a connection on $E$, determines a superconnection
on $\pi^*E$ \cite{Q}. Quillen's superconnection $\nabla+L$ encodes
all the information necessary to define the Chern character in a
single object. We shall
denote it by $\nabla_L$.

The supertrace of $\exp\nabla_L^2$ is a mixed differential form
which enjoys the major properties of the ordinary Chern
character: it is closed and its cohomology class depends
only on the underlying vector bundle, not on $\nabla$ or $L$.
But rather than passing to its cohomology class on $T^*M$,
we take advantage of the fact that this form is rapidly decreasing on
the fibers of $T^*M$.  We get this convenient property
by sacrificing the traditional ${\frac {-1}{2\pi i}}$
factor, an error which we shall also  address.
Thus, we define a current on $\Omega^*M$ by the Poincare Duality formula
$$PD\co\eta \mapsto \int_{T^*M}\pi^*(\eta)\operatorname{tr_s}\exp\nabla_L^2.$$
If we expand this dual  ``Chern character current''
in the Taylor series, we obtain:

$$PD\left[\operatorname{tr_s}\exp\nabla_L^2\right]
   =PD\left[\operatorname{tr_s}(1+\nabla_L+{\frac 1 {2!}}
\nabla_L^2+\ldots)\right]$$
The terms on the right are closed forms and their cohomology
classes depend on the isomorphism type of $\pi^*E$ only.
However, they are not rapidly decreasing and we cannot
form the dual currents by coupling them with the pullback of
an arbitrary smooth form on $M$.

We resolve this issue through {\it analytic regularization}, to be
addressed in section \ref{section:Form}, and thus obtain a new
formula for the dual Chern character. Briefly, our main result can
be stated as follows:

{\bf Theorem. } For $R>0$, let $Y_R$ be the $R$-tubular
neighborhood of
the zero section of $T^*M$. Then:
\begin{align*}
\int_{T^*M}\pi^*(\eta)&\operatorname{tr_s}\exp\nabla_L^{2}\\&=
    \lim_{R\to 0}
\sum_{z\in{\C}}Res|_z\Big[\Gamma(z)
\int_{T^*M\backslash Y_R}\pi^*(\eta)\operatorname{tr_s}(-\nabla_L)^{-2z}\Big].
\end{align*}
The limit arises due to divergence of negative powers of $L$ near zero.


We now proceed to briefly describe the Connes-Moscovici
local Chern character, i.e., the noncommutative Chern character.

The Chern character in \cite{CM} is a periodic cyclic cohomology
class in $HPC^*$. This cohomology is constructed from spaces of
multilinear functionals on $A$. See \cite{GVF}, chapter 10, for
construction of $HPC^*$. The Connes-Moscovici Chern character is
represented by a sequence of multilinear functionals on $A$ which
we proceed to describe. The n-th term of the sequence is  zero for
$n=1,3,5,\ldots$; for $n=0,2,4\ldots$, let $a_0, a_1,\ldots a_n$
be the elements of $A$. Let $\operatorname{Tr_s}$ be  the
supertrace and let $k$ be the variable running through {\it all}
$n$-multiindices with nonnegative integer entries.
\begin{align}\label{section:cmchfmla}
\phi(a_0,a_1,\ldots a_n)
=&\sum_kC_{nk} Res|_{z=0}\Gamma(z+|k|+{\frac n2})\nonumber\\&
\,\,{\scriptstyle \times}
\operatorname{Tr_s}
\left({\scriptstyle a_0}\delta^{(k_1)}
({\scriptstyle[D,a_1]})
\ldots \delta^{(k_n)}({\scriptstyle[D,a_n]})
D^{-2(z\!+\!|k|\!+\!{\frac n2})}\right),\nonumber
\end{align}
where
$$C_{nk}={\frac {(-1)^{|k|}\Gamma({\frac n2}+|k|)}
        {k!(k_1+1)(k_1+k_2+2)\ldots(k_1+k_2+\ldots k_n+n)}}.$$

Also, note that the trace $Tr_s(\ldots\,D^{-2(z+|k|+{\frac n2})})$
 {\it must be replaced by its meromorphic continuation in $z$, before
we take the residues}.
Existence of such a continuation is also implied by certain axioms and is
not trivial at all. Indeed, as it stands,
the operators whose trace we are taking  typically
fail to be bounded, let alone trace class.

To us, the most important fact about the Connes-Moscovici
formula is that this Chern character is a
sum of residues of $Tr_s(\ldots\,D^{-2(z+|k|+{\frac n2})})$
  times the gamma function. Our main
result expresses the Quillen's representative of the (classical)
Chern character as a sum of very similar quantities.
Much like the proof in \cite{CM}, our argument hinges on the Mellin
transform.

\section{Superconnections and Chern Classes.}\label{section:sch}

In this section, we give an overview of superconnections according
to Quillen \cite{Q}. This notion shall be used to define Chern character in
the spirit of Chern-Weil theory. Let $M$ be a smooth manifold with
no boundary. Let $E$ be a smooth complex ${\Z}_2$-graded vector bundle over $M$.
We work with the vector bundle $\Lambda^*M\otimes E.$

\begin {defn}(Quillen, \cite{Q})
Let ${\alg}(E)$
be the space of smooth sections of the vector bundle
$\Lambda^*M\otimes End(E).$
It naturally inherits the  ${\Z}_2$-grading from the fibers.
\end{defn}

\begin{defn}\label{section:grmult}
Let $\omega,\,\nu$ be homogeneous differential forms
on $M$, let $T, S\in \Gamma^\infty End(E)$ be homogeneous
(purely even or purely odd) endomorphisms of $E$ and let
$s\in\Gamma^\infty(\Lambda^*M\otimes E)$
be a homogeneous section.
We define the {\it graded multiplication} on ${\alg}(E)$
by:
\begin{equation}\label{section:oml1}
(\omega\otimes T) \,(\nu\otimes S)=_{def}(-1)^{deg(\nu)deg(T)}\omega\nu
  \otimes TS
\end{equation}
Also, the action of ${\alg}(E)$ on $\Gamma^\infty(\Lambda^*M\otimes E)$
is defined by:
\begin{equation}\label{section:oml2}
(\omega\otimes T )\,(\nu\otimes s)=_{def}
(-1)^{deg(\nu)deg(T)}T\nu\otimes Ts
\end{equation}
\end{defn}

\begin{lemma}
The equation (\ref{section:oml1}) makes ${\alg}(E)$ an
associative  superalgebra. Also, (\ref{section:oml2})
defines an algebra action of ${\alg}(E)$ on
$\Gamma^\infty(\Lambda^*M\otimes E).$ In fact, this makes ${\alg}(E)$
a subalgebra of $\Gamma^\infty End(\Lambda^*M\otimes E)$
in the sense that no nonzero element of ${\alg}(E)$ kills
everything.
\qed
\end{lemma}

The condition (\ref{section:oml2}) is called $\Omega^*$-{\it linearity.}
It turns out that  $\Omega^*$-linear endomorphisms of
$\Lambda^*M\otimes E$ are precisely the elements of ${\alg}(E)$.

\begin{lemma} The algebra
${\alg}(E)$ is the ${\C}$-span of homogeneous
smooth sections of $End^{\pm}(\big(\Lambda^*M\big)\otimes E)$
which are in addition $\Omega^*$-linear.
\qed
\end{lemma}
\begin{remark}
{\rm Observe that ${\alg}(E)$  is naturally a left $\Omega^*M$-module.}
\end{remark}
\begin{remark}{\rm
The definition \ref{section:grmult} resembles the usual definition
of multiplication on a tensor product of two algebras. However, it
takes the ${\Z}_2$-grading into account.  Such products of
superalgebras are called {\it graded tensor products} and are
denoted by $\grtimes$. }
\end{remark}

Next, observe that the ${\Z}_2$-grading on $E$ naturally induces
one on the space of sections $\Gamma^\infty(\Lambda^*M\otimes E)$,
which makes it possible to speak of even and odd linear
endomorphisms of this space. The even ones preserve the ${\Z}_2$
grading and the odd ones switch it.  We are now ready to define
superconnections.

\begin{defn} Let $E$ be a smooth ${\Z}_2$-graded complex vector bundle.
In the spirit of \cite{Q}, we define a {\it
superconnection on} $E$ to be an odd ${\C}$-linear map
$$\nabla\co\Gamma^\infty(\Lambda^*M\otimes E)
     \to \Gamma^\infty(\Lambda^*M\otimes E),$$
which satisfies the so-called
{\it graded Leibniz rule}:
$$\nabla(\omega s)=(d\omega)s+(-1)^{deg (\omega)}\omega\nabla s.$$
The {\it curvature of } $\nabla$ is defined as $\nabla^2.$
\end{defn}

\begin{lemma}
The curvature is a globally well-defined,
even element of ${\alg}(E).$\qed
\end{lemma}

\begin{lemma}
Locally, any superconnection is always of the form \\$d\otimes id +
\theta$, where $d$ is the de Rham differential and $\theta$ is a
local odd section of $\Lambda^*M\otimes End(E).$ \qed
\end{lemma}

\begin{corollary}
Given a superconnection $\nabla$ and an  arbitrary odd
element $L$ of ${\alg}(E)$, $\nabla+L$ is
also a superconnection.
\end{corollary}

\begin{lemma}Any connection which respects the ${\Z}_2$-grading
is a superconnection.
\end{lemma}
{\bf Proof:} Such a connection can be locally written as
$d+\theta$, where $\theta=(dx_i\Gamma^k_{ij})_{jk}$ is a matrix of
1-forms. Since the connection respects the grading, the matrices
$(\Gamma_{ij}^k)_{ik}$ define even endomorphisms of $E$. Presence
of 1-forms, therefore, makes $\theta$ odd.\qed

\begin{defn}
We denote the superconnection $\nabla+L$  by $\nabla_L$.
\end{defn}

\begin{defn}
Suppose $W$ is a finite-dimensional ${\Z}_2$-graded
complex vector space and
$A=End(W)$. We define the {\it supertrace} $\operatorname{tr_s}^{\C}\co A\to {\C}$
by the following equation:

$$\forall  f=\left(
       \begin{array}{cc}
       f_{11}&f_{12}\\ f_{21}&f_{22} \\
      \end{array} \right) \in A,\,\,
\operatorname{tr_s}^{\C}(f)=_{def} \operatorname{tr}(f_{11})-
\operatorname{tr}(f_{22}),$$
where the traces on the right hand side are the usual traces
of operators from $W^+$ and $W^-$ into themselves.
Let $R$ be another ${\Z}_2$-graded algebra.
A typical simple tensor in $R\grtimes A$ has the form
$$ f=\left(
       \begin{array}{cc}
       r_{11}f_{11}&r_{12}f_{12}\\ r_{21}f_{21}&r_{22}f_{22} \\
      \end{array} \right).$$
On simple tensors, we define $\operatorname{tr_s}^{R}\co
R\grtimes A\to {R}$ by
$$r_{11}{\operatorname tr}(f_{11})-(-1)^{deg(r_{22})}
  r_{22}{\operatorname tr}(f_{22}).$$
\end{defn}
If $W$ is a fiber of some vector bundle $E$,
this definition extends naturally to sections.

We proceed toward the definition of Chern character. To that end,
we need the graded version of the trace. Being ${\Z}_2$-graded,
$End(E)$ has a supertrace $\operatorname{tr}^{\C}_s\co End(E)\to{\C}$ which means
that on $\Lambda^*M\grtimes End(E)$ the $\Lambda^*V$-valued
supertrace $\operatorname{tr}^{\Lambda^*M}\co\Lambda^*M\otimes (E)\to
\Lambda^*M $ makes sense.

\begin{defn}\label{theorem:gchdef}
The $2k$-th component of the {\it Chern character form} is
defined as the following differential form:

$$ch_{2k}(E)={\frac 1 {k!}}
 \operatorname{tr_s}^{\Lambda^*M}(\nabla^{2k})\in \Omega^{2k}M.$$

The {\it total Chern character form} is
$$ch(E)=\sum_{k=0}^\infty ch_k(E).$$
 We usually write this
form as $\operatorname{tr_s}^{\Lambda^*M}(\exp(\nabla^2))$.
\end{defn}

\begin{theorem}(Quillen, \cite{Q}).
The  series defining $ch(E)$ converges.\qed
\end{theorem}

\begin{theorem}\label{theorem:gchthm} (\cite{Q}) $ch_k$ is closed
and its cohomology class is
independent on a choice of the superconnection.
This class is an  invariant of isomorphism classes
of complex vector bundles.\qed
\end{theorem}

\section{Formulation of the Problem.}\label{section:Form}
Here we state the main theorems of the paper. The proofs
shall be given in subsequent sections.

The general geometric setup for our problem is the following. Let
$M$ be a compact, $n$-dimensional, smooth manifold. Let
$E\rightarrow M$ be a smooth, ${\Z}_2$-graded complex vector bundle.
We assume that $E$ and $TM$ are provided with metrics.
Let $D$ be an elliptic, odd, first-order, selfadjoint differential
operator on $E$. Finally, let $\pi\co T^*M\to M$ be the
standard projection map of the cotangent bundle,
and let $\nabla$ be a superconnection on $\pi^*E$
which arises as a pullback of some superconnection
$\nabla'$ from $E$.

Our result is motivated by the Connes-Moscovici
formula. We express the right-hand side of the Atiyah-Singer
Index Theorem in a way which resembles their
residue cocycle. Namely, we shall sum over all the
residues of the expression:
$$\Gamma(z)\int_{T^*M} \pi^*(\eta)\nabla_L^{-2z}.$$

The Index Theorem using Quillen's Chern charactrer can be stated as:
$$Ind(D) = \sum_\kappa\big({\frac {-1}{2\pi i}}\big)^{2n-\kappa}
\int_{T^*M}\pi^*(Todd[TM\otimes {\C}])_\kappa ch(L).$$
Here we need to prove that the integral
(lemma \ref{section:expdecr}) converges and that $ch(L)$ may indeed be used
in place of the ordinary chern character (theorem \ref{theorem:0}).
The factor of $({\frac {-1} {2\pi i}})^{2n-\kappa}$ is to correct for
the error introduced by leaving the $2\pi i$ out of Quillen's
definition of Chern chracter.

The right-hand side is really concerned with the Poincare Dual
of $ch(L)$, i.e. with the following  linear functional
$$PD[ch(L)]\co\eta\mapsto\int_{T^*M}\pi^*(\eta)ch(L).$$

\begin{lemma}\label{section:expdecr}
The quantity $\|ch(L)\|=\|\operatorname{tr_s}\exp\nabla_L^2\|$
decays exponentially along the fibers of $T^*M$.
In fact, there are positive constants $C$ and $K$ such that
$$\|ch(L)\|\le Ce^{-K\rho^2},$$
where $\rho$ is the radial coordinate on the fibers
of $T^*M$ obtained from the Riemanian metric.
\end{lemma}

\begin{corollary} Since $M$ is compact, the integral
$\int_{T^*M}\pi^*(\eta)ch(L)$ converges.
\end{corollary}

The necessary estimates for this lemma are provided in
section 3 of \cite{Q}. Essentially, it is true because
$D$ is selfadjoint, so that $L$ is antiselfadjoint
and $L^2$ is negative-definite, which is where the ellipticity of $D$ comes
in. Also, $L^2$ increases
polynomially on the fibers of $T^*M$. When we exponentiate
$\nabla^2_L=\nabla^2+[\nabla,L]+L^2$, the resulting expression
decreases as $e^{-Kr^2}$.

Next consider  the Taylor expansion of  $ch(L)$:
$$\operatorname{tr_s}\exp\nabla_L^2=
  \operatorname{tr_s}(1+\nabla_L+{\frac 1 {2!}}
\nabla_L^2+\ldots).$$

Although the Poincare dual of the left-hand side converges,
the duals the individual Taylor terms do not,
due to polynomial increase of $L$ along the fibers
of $T^*M$.

We get around this difficulty by analytic regularization.
The general idea is that if we replace the integer power $k$
of $\nabla^2_L$ by a complex number $-z$, where $Re(z)\gg 1$,
we also can replace the integral
${\frac 1{k!}}\int_{T^*M} \pi^*(\eta)\nabla^{2k}_L$
with the following expression:
\begin{equation}\label{section:whydiverge}
Res|_{z=-k}\Gamma(z)\int_{T^*M} \pi^*(\eta)\nabla^{-2z}_L.
\end{equation}

Before we can take this residue, though, we need to
pass to the meromorphic extension of
\begin{equation}\label{section:nablz}
\int_{T^*M} \pi^*(\eta)\nabla^{-2z}_L.
\end{equation}
In particular, we need to prove that such an
extension exists (theorem \ref{theorem:1}).

However, there is yet another difficulty. Let
$|\,\,|$ denote the fiberwise norm on $T^*M$.
Define:
$$Y_R=_{def}\{\xi\in T^*M\co|\xi|<R\}.$$
$$X_R=_{def}Y_R^c.$$

The integral in (\ref{section:nablz})
would not really converge, if taken over all of $T^*M$,
because $L^{2}$ is a symbol of a second-order
differential operator. Hence, $L^2$ is a homogeneous quadratic polynomial
in the vertical coordinates $\xi$ of the cotangent
bundle. It therefore vanishes at the zero section on $T^*M$.

The problem with divergence at infinity
shall be resolved by taking $Re(z)>0$, meromorphically extending the
integral to the whole complex plane and taking residues.

We deal with divergence at $Y_R$, in the following way.
First, we shall replace the integral over $T^*M$ by that over $X_R$.
Second,  we shall take  the residues in $z$. Third, we shall take the
limit as $R\to 0$. That is,
we use:
$$\lim_{R\to 0}\sum_{z\in{\C}}Res|_z\Gamma(z)
\int_{X_R}\pi^*(\eta)\operatorname{tr_s}(-\nabla_L^2)^{-z}.$$
It turns out that this quantity is well-defined and is equal to the
value of the current $PD(\operatorname{tr_s}\exp\nabla_L^2)$ on $\eta$.

We now formally state the main results
of the present work. First, we summarize the hypotheses
which apply in all the theorems in the sequel:
{\it Let $M$ be a compact, smooth $n$-manifold with no boundary.
Let $\pi\co T^*M\to M$ be the cotangent bundle and let
$E\to M$ be a ${\Z}_2$-graded smooth vector bundle.
Suppose also that $D$ is an odd, elliptic, first order
selfadjoint differential operator on $E$.
Thus, the symbol $L$ of $D$ is an odd endomorphism of $\pi^*E$
which is invertible everywhere but at the zero section $M\subseteq T^*M$
and pointwise anti-selfadjoint.
Let $\nabla$ be the pullback of some connection $\nabla'$
on $E$ via $\pi$. Assume also that $\nabla'$, and hence $\nabla$,
respects the ${\Z}_2$-grading.}

\begin{theorem}\label{theorem:0}
 Let $\eta$ be a closed, smooth differential form on $M$.
 Consider the following
integral:
$$I(\eta)=\int_{T^*M}\operatorname{tr_s}
\pi^*(\eta)\exp{\nabla_L^{2}}.$$
\begin{itemize}
\item[a)] It vanishes if $\eta$ is exact.
\item[b)] It is independent of the particular choice of $\nabla'$.
\item[c)] Let $\beta=-\int_0^\infty\operatorname{tr_s}
           \exp\nabla_{tL}^2Ldt.$
          For any $R>0$, the following holds
          on the interior of $X_R$:
   $$d\beta=\operatorname{tr_s}\exp\nabla^2.$$
          Thus, the pair $(\operatorname{tr_s}\exp\nabla^2,\beta)$
          defines a relative cohomology class in
   $\scriptstyle H^*(T^*M,int\,X_R)$.
\item[d)]
$I(\eta)=\int_{_{T^*M\backslash X_R}}
\pi^*(\eta)\operatorname{tr_s}\exp\nabla^2\,\,-\,\,
          \int_{_{\del(T^*M\backslash X_R)}}\pi^*(\eta)\beta.$\\
Hence, $I(\eta)$ yields the same result as pairing of $\eta$
 with the relative cohomology class defined by
$(\operatorname{tr_s}\exp\nabla^2,\beta).$

\end{itemize}
\end{theorem}

Observe that $\pi^*\co H^*M\cong H^*(T^*M)$ and that
$\operatorname{tr_s}\exp\nabla^2$ and $\pi^*(chE^+-chE^-)$ are
cohomologous. That is, $\operatorname{tr_s}\exp\nabla^2$
determines the same difference Chern class.

\begin{theorem} \label{theorem:1}
For any $R>0$ and for any
$\eta\in\Omega^*M$ and for any
$z\in{\C}$ with $Re(z)\gg0$, the following integral converges:
\begin{align}
\int_{X_R}\pi^*(\eta)\operatorname{tr_s}(-\nabla_L^2)^{-z}.
\end{align}
Further, it has a meromorphic extension to ${\C}$ of the form
$$\sum_K{\frac {R^{K+1-2z}}{K+1-2z}}A_K,$$
where $A_K$ are constants.
\end{theorem}

In fact, we do have a classification of the poles. There are
finitely many of them and they are located at negative integers
or half-integers.

\begin{lemma} \label{theorem:poles}
Let $\eta\in\Omega^{\kappa}T^*M$. Then the integral
$$\int_{X_R}\pi^*(\eta)\operatorname{tr_s}(-\nabla_L^2)^{-2z}$$
can only have a nonzero residue at the point
$(\kappa-2n)/2$. Further, if $n$ is even and the residue
is nonzero, then $(\kappa-2n)/2$ must be an integer. If $n$
is odd,  $(\kappa-2n)/2$ must be a half-integer. In either case,
nonzero residues occur only for even $\kappa$.
\end{lemma}

This lemma says that at each particular point  $z\in {\C}$,
the residue is a homogeneous current. That is, it vanishes
on all the forms except possibly for those of some given degree.

\begin{theorem} \label{theorem:2}
Let $\nabla$ be a superconnection which has been pulled back
from $E$ via $\pi$. For any  $\eta\in\Omega^*M$
\begin{multline}\label{theorem:21}
\int_{T^*M}\!\!\!\pi^*(\eta)\operatorname{tr_s}\exp{\nabla_L^{2}} =
\sum_{z\in{\C}}\lim_{R\to 0}Res|_z\Gamma(z)\!\!
\int_{X_R}\pi^*(\eta)\operatorname{tr_s}
({ -\nabla_L^2})^{-z}.
\end{multline}

All but finitely many residues on the
right-hand side vanish as $R\to 0$.
\end{theorem} \qed

\begin{corollary}
For each $z$, the following defines a closed current on $\Omega^*M$:
$$R_z\co \eta\mapsto\lim_{R\to 0}Res|_z\Gamma(z)
\int_{X_R}\pi^*(\eta)\operatorname{tr_s}
({ -\nabla_L^2})^{-z}.$$
\end{corollary}
{\bf Proof:}
By theorem \ref{theorem:gchthm},
$\operatorname{tr_s}\exp\nabla^2_{tL}$ is a closed form. It is
also rapidly decreasing on the fibers of $T^*M$, so
for an exact form $\eta=d\omega$ on $M$,
Stokes theorem yields:
$$\int_{T^*M}\pi^*d\omega\operatorname{tr_s}\exp\nabla^2_{tL}=
\pm \int_{T^*M}\pi^*\omega d\operatorname{tr_s}\exp\nabla^2_{tL}=0.$$
The rapid decay property assures that there is no boundary term.
Thus, $\operatorname{tr_s}\exp\nabla^2_{tL}$ induces a
closed current on $\Omega^*M$. By lemma \ref{theorem:poles}, for
each  $\kappa$, $R_z$ either vanishes on
or agrees with the current induced by
$\operatorname{tr_s}\exp\nabla^2_{tL}$. In either case,
$R_z$ is a closed current on $\Omega^\kappa M$.
\qed\\

Just as the computation in \cite{CM}, our
proof of theorem \ref{theorem:2} hinges on Mellin
transform. The simplest example of a Mellin transform is the well-known formula,
valid for $Re(\sigma)>0$:
$$\int_0^\infty e^{-\sigma t}t^{z-1}dt=\sigma^{-z}\Gamma(z).$$
It says that $\sigma^{-z}\Gamma(z)$  is the {\it Mellin transform
} of $e^{-\sigma t}$. (See \cite{Arfken} for details).
Connes and Moscovici apply the same
transform to the so-called {\it JLO cocycles} \cite{Q3} in order
to obtain the residue cocycle. Let $\Sigma_k$ be the standard
$k$-simplex in ${\R}^{k+1}$ with coordinates
$u_0,u_1\ldots\,u_{k-1}$.  Using the notation of section
\ref{section:intro}, the JLO cocycles are comprised of multilinear
functionals on a *-algebra $A$ given by
$$\psi_{JLO}^t(a_0, a_1,\ldots, a_k)=
\operatorname{Tr_s}\int_{\Sigma_k}a_0e^{-u_0tD^2}a_1e^{-u_1tD^2}\ldots
a_{k-1}e^{-u_{k-1}tD^2}a_kd{\mathbf u}.$$
We, however, apply Mellin transform to
Quillen's Chern character, which is an exponential,
and obtain complex powers of the curvature.

\section{Proof of Theorem \ref{theorem:0}.}\label{section:pt0}

{\bf Theorem \ref{theorem:0}.} {\it  Let $\eta$ be a closed
differential form on $M$. Under the hypotheses outlined in section
3, consider the following integral:
$$I(\eta)=\int_{T^*M}\operatorname{tr_s}\pi^*(\eta)\exp{\nabla_L^{2}}.$$
\begin{itemize}
\item[a)] It vanishes if $\eta$ is exact. \item[b)] It is
independent of the particular choice of $\nabla'$.
\item[c)]Let
$\beta=-\int_0^\infty(\exp\nabla_{tL}^2)Ldt$.
           For any $R>0$, the following equation holds
           on the interior of $X_R$:
           $$d\beta=\operatorname{tr_s}\exp\nabla^2.$$
          Thus, the pair $(\operatorname{tr_s}\exp\nabla^2,\beta)$
          defines a relative cohomology class in
   $\scriptstyle H^*(T^*M,int\,X_R)$.
\item[d)]
$I(\eta)=\int_{_{T^*M\backslash X_R}}\pi^*(\eta)\operatorname{tr_s}\exp\nabla^2\,\,-\,\,
          \int_{_{\del(T^*M\backslash X_R)}}\pi^*(\eta)\beta.$\\
Hence, $I(\eta)$ yields the same result as pairing of $\eta$
 with the relative cohomology class defined by
$(\operatorname{tr_s}\exp\nabla^2,\beta).$

\end{itemize}
}

The content of this theorem is really due to \cite{Q}.  For part
(a), assuming $\eta=d\omega$, we compute:
$$\int_{_{Y_R}}\pi^*(d\omega)\operatorname{tr_s}
       \exp{\nabla_L^{2}}=
\int_{_{\partial Y_R}}\pi^*(\omega)\operatorname{tr_s}
       \exp{\nabla_L^{2}}.
$$  The right-hand side vanishes as $R\to\infty$.

For (b), the fact that the cohomology class of
$\operatorname{tr_s}\exp\nabla_{L}^2$
 is independent on  $\nabla'$ or $L$ is not enough.
We need to prove that if we replace the connection $\nabla'$ on
$E$ with some ${\tilde \nabla}'$ , so that on $\pi^*E$ we have
$\tilde\nabla=\pi^*\tilde\nabla'$, then  there exists a
differential form $\beta_1$, rapidly decreasing on the fibers of
$T^*M$ and such that
 $$d\beta_1=\operatorname{tr_s}\exp\nabla_{L}^2
  -\operatorname{tr_s}\exp\tilde\nabla_L^2.$$
We present the so-called ''homotopy'' argument.

Suppose first there is some connection $\nabla_t$ on
$\pi^*E$ which depends
on a parameter $t$ in a differentiable way.
The example we have in mind is:
$$\nabla_t=\pi^*(t\nabla'+(1-t){\tilde \nabla}')+L=
t\nabla+(1-t){\tilde \nabla}+L.$$

Here, $t$ is a coordinate
on the manifold $T^*M\times {\R}$. We take the pullback vector
bundle $F=pr_1^*(\pi^*E)$, on $T^*M\times {\R}$ and the
following  defines a
connection $D_t$ on $F$: $$D_t=_{_{def}} \nabla_t+dt\del_t.$$

Then $D_t^2=\nabla_t^2+dt\dot{\nabla_t}$ and more generally:
\begin{align}\label{section:dt1}
D_t^{2k}=\nabla_t^{2k}+
\sum_{j=0}^{k-1}\nabla_t^{2j}dt\dot{\nabla_t}
      \nabla_t^{2(k-j-1)}=\mu_k+dt\,\,\nu_k,
\end{align}
where $\mu_k$ and $\nu_k$ are unambiguously defined
by the above equation.

Let $d'$ denote the de Rham differential on $T^*M\times{\R}$
and denote the one on $T^*M$ by simply $d$.
This way, $d'=d+dt\del_t$.
By theorem \ref{theorem:gchthm}, $\operatorname{tr_s}D_t^{2k}$ is closed:
$$d'\operatorname{tr_s}D_t^{2k}=d'\mu+d(dt\,\,\nu)=0$$
so that $\del_t\mu_k=d\nu_k$ and $\mu_k|_{t=1}-\mu_k|_{t=0}=d\int_0^1\nu dt$,
which means that:
$$\nabla_1^{2k}-\nabla_0^{2k}=d\int_0^1
\sum_{j=0}^{k-1}\nabla_t^{2j}\dot{\nabla_t}\nabla_t^{2(k-j-1)}dt.$$
Or, taking supertraces and keeping in mind that supertraces
kill supercommutators:
$$\operatorname{tr_s}(\nabla_1^{2k}-\nabla_0^{2k})=
d\operatorname{tr_s}\int_0^1
k\nabla_t^{2k}\dot{\nabla_t}dt.$$

This implies that
\begin{align}\label{section:0partb0}
\operatorname{tr_s}\exp \nabla_1^2-
\operatorname{tr_s}\exp \nabla_0^2=
d\int_0^1\operatorname{tr_s}\exp \nabla_t^2\dot{\nabla_t}dt.
\end{align}
This equation applies
to  any $\nabla_t$ which depends on $t$ differentiably.

In our particular example, $$\nabla_t=
\pi^*(t\nabla'+(1-t){\tilde \nabla}')+L,$$
so that $\nabla_0=\tilde\nabla_L$ and $\nabla_1=\nabla_L$.
 It follows that
\begin{align}\label{section:0partb1}
\operatorname{tr_s}\exp \nabla^2_L-
\operatorname{tr_s}\exp \tilde\nabla^2_L=
d\int_0^1\operatorname{tr_s}\exp \nabla_t^2\dot{\nabla_t}dt,
\end{align}  and we may define $\beta_1$ as
$\int_0^1\operatorname{tr_s}\exp\nabla_t^2\dot{\nabla_t}dt.$\\
Because $\exp \nabla_t^2$ is exponentially decreasing along the
fibers of $T^*M$ for each fixed $t$, the rapid decay property of
$\beta_1$ is easy to prove.

To proceed with (c), we change our definition of $\nabla_t$:
$$\nabla_t=\nabla+tL,$$ and $D_t=\nabla_t+dt\del_t$
on $T^*M\times {\R}.$
This does not affect the
fact that $\del_t\mu_k=d\nu_k$ and we have:
 $$\del_t\nabla_t^{2k}=d\sum_{j=0}^{k-1}\nabla_t^{2j}
                 \dot{\nabla_t}\nabla_t^{2{k-1}},$$
or, taking supertraces:
\begin{align}
\operatorname{tr_s}\del_t\nabla_t^{2k}&=
    d\operatorname{tr_s}k\nabla_t^{2(k-1)}
                 \dot{\nabla_t}
\end{align}

Just as in part (b), we obtain:
\begin{align}\label{section:delt}
\del_t\operatorname{tr_s}\exp\nabla_t^{2}=
    d\operatorname{tr_s}\exp \nabla_t^{2}L.
\end{align}

But as long as $L$ is invertible, $\|\exp \nabla_t^{2}\|\to 0$ as
$t\to 0$ (since, there are non-zeroes among the eigenvalues of
$L^2$). Hence, (integrating (\ref{section:delt})) yields:
\begin{align}\label{section:0partc1}
\operatorname{tr_s}\exp\nabla_L^{2}&=
\operatorname{tr_s}\exp\nabla_1^{2}\\
    &=-d\int_1^\infty\operatorname{tr_s}\exp \nabla_t^{2}
                 \dot{\nabla_t}dt.
\end{align}
So, the above equation holds on $X_R$. Keeping in
mind (\ref{section:0partb0}), which can be restated here as:
\begin{align}\label{section:0partc2}
\operatorname{tr_s}\exp \nabla^2_L-
\operatorname{tr_s}\exp \nabla^2=
d\int_0^1\operatorname{tr_s}\exp\nabla_{tL}^2Ldt,
\end{align}
we see that on $X_R$:
\begin{align}\label{section:0partc3}
\operatorname{tr_s}\exp \nabla^2=
-d\int_0^\infty\operatorname{tr_s}\exp\nabla_{tL}^2Ldt,
\end{align}

For (d), recall how does one integrate compactly supported cohomology
classes defined by pairs.
If we have a pair $(\operatorname{tr_s}\eta\exp \nabla^2,\eta\beta)$
as above, $\beta$ being equal to
$\int_0^\infty\operatorname{tr_s}\exp\nabla_{tL}^2Ldt$,
and $\eta$ being closed, then:
\begin{align}
\langle[T^*M],[(\operatorname{tr_s}\eta\exp \nabla^2,\eta\beta)]\rangle=_{def}
\int_{Y_R}\eta\operatorname{tr_s}\exp\nabla^2 -\int_{\del
Y_R}\eta\beta.
\end{align}

Similarly, for
$\beta_2=-\int_1^\infty\operatorname{tr_s}\exp\nabla_{tL}^2Ldt$ we
already have:
\begin{align}
\langle[T^*M],[(\operatorname{tr_s}\eta\exp \nabla_L^2,\eta\beta_2)]\rangle=_{def}
\int_{Y_R}\eta\operatorname{tr_s}\exp\nabla_L^2 -\int_{\del
Y_R}\eta\beta_2.
\end{align}

Combining (\ref{section:0partc1}) and (\ref{section:0partc2}) we see
that these quantities are equal. Taking $R\to\infty$, due to
exponential decay, $\int_{\del
Y_R}\eta\beta\to 0$. It follows that
$$\int_{T^*M}\eta\operatorname{tr_s}\exp\nabla_L^2=
\int_{Y_R}\eta\operatorname{tr_s}\exp\nabla^2 -\int_{\del
Y_R}\eta\beta.$$

\section{Proof of Theorem \ref{theorem:1}.}\label{section:pt1}
{\bf Theorem \ref{theorem:1}}\begin{it}
Under the hypotheses outlined in section \ref{section:Form},
for any $R>0$ and for any
$\eta\in\Omega^*M$ and for any
$z\in{\C}$ with $Re(z)\gg0$, the following integral converges:
\begin{align}\label{theorem:1eq}
\int_{X_R}\pi^*(\eta)\operatorname{tr_s}(-\nabla_L^2)^{-z}.
\end{align}
Further, it has a meromorphic extension to ${\C}$ of the form
$$\sum_K{\frac {R^{K+1-2z}}{K+1-2z}}A_K,$$
where $A_K$ are constants.
\end{it}

First, we set up some notation. Over each coordinate chart
$U_\alpha$ of $M$, we may define a pullback chart
$V_\alpha=\pi^{-1}(U_\alpha)$ of $T^*M$. It has horizontal
coordinates $x=(x^1,\ldots\,,x^n)$, which are just the coordinates
of $M$, and vertical ones $\xi=(\xi^1,\ldots\,,\xi^n)$. On the
fibers of $T^*M$, let $\rho$ and  $\Xi$ be the spherical
coordinates. We can take $\rho=|\xi|$ to be the coordinate induced
by the metric $g$. Let $S_\rho^*M$ denote the sphere bundle of
$T^*M$ of radius $\rho$. The theorem follows by direct computation
from the following proposition.

\begin{proposition}\label{section:details}
Over each chart $X_R\bigcap V_\alpha$, there exist smooth local
sections $\Theta^K(z,x,\rho,\Xi)$ of
$\Lambda^*T^*M\grtimes End(\pi^*E)$ whose coordinate expressions,
in fact, do not
depend on $\rho$ (i.e., they are pullbacks
from the unit sphere bundle via the obvious map  $X_R\to S^*M$),
so we write $\Theta^K(z,x,\Xi)$.
For all $z$ with $Re(z)\gg0$, for all
$R>0$ and $\eta\in\Omega^*M$, these sections  $\Theta^K(z,x,\Xi)$
satisfy:
\begin{align}\label{section:details00}
\int_{_{X_R\bigcap V_\alpha}}
\!\!\!\!\!\!\!\!\!\!\!\!\!\!\pi^*(\eta)\operatorname{tr_s}(-\nabla_L^2)^{-z}\!&=\\
\nonumber
&\!\int_R^\infty\!\!\!\!\rho^{-2z+K}d\rho\!\!\int_{_{S_\rho^*M\bigcap
V_\alpha}} \!\!\!\!\!\!\!\!\!\!\!\!\pi^*(\eta)
  \sum_{K}\operatorname{tr_s}
\Theta^K(z,x,\Xi) d\Xi.
\end{align}
The sum in the right-hand side is finite.
Further, each $\Theta^K(z,x,\Xi)$ extends to an entire function in $z$
which is also smooth in $x$ and $\Xi$.
\end{proposition}

This proposition follows, essentially, by separation of powers of
$\rho$ in the coordinate expression of $\nabla_L^2$.

 {\bf Proof of
proposition} \ref{section:details}. The major steps in the proof
are the following. First, for a suitable contour $\gamma$, we
express $(-\nabla_L^2(p))^{-z}$ at each $p\in X_R$
by holomorphic functional calculus:
\begin{align}\label{section:details0}
(-\nabla_L^2)^{-z}={\frac 1{2\pi i}}\int_\gamma
\lambda^{-z}(\lambda +\nabla^2_L)^{-1}d\lambda.
\end{align}
Convergence of the integral is obvious and the fact that $\gamma$
may be used instead of the usual counter-clockwise oriented curve
follows by standard argument as in figure 1.
Note that $\gamma$ does not depend on $p$.

\begin{figure}
\centerline {\psfig {figure=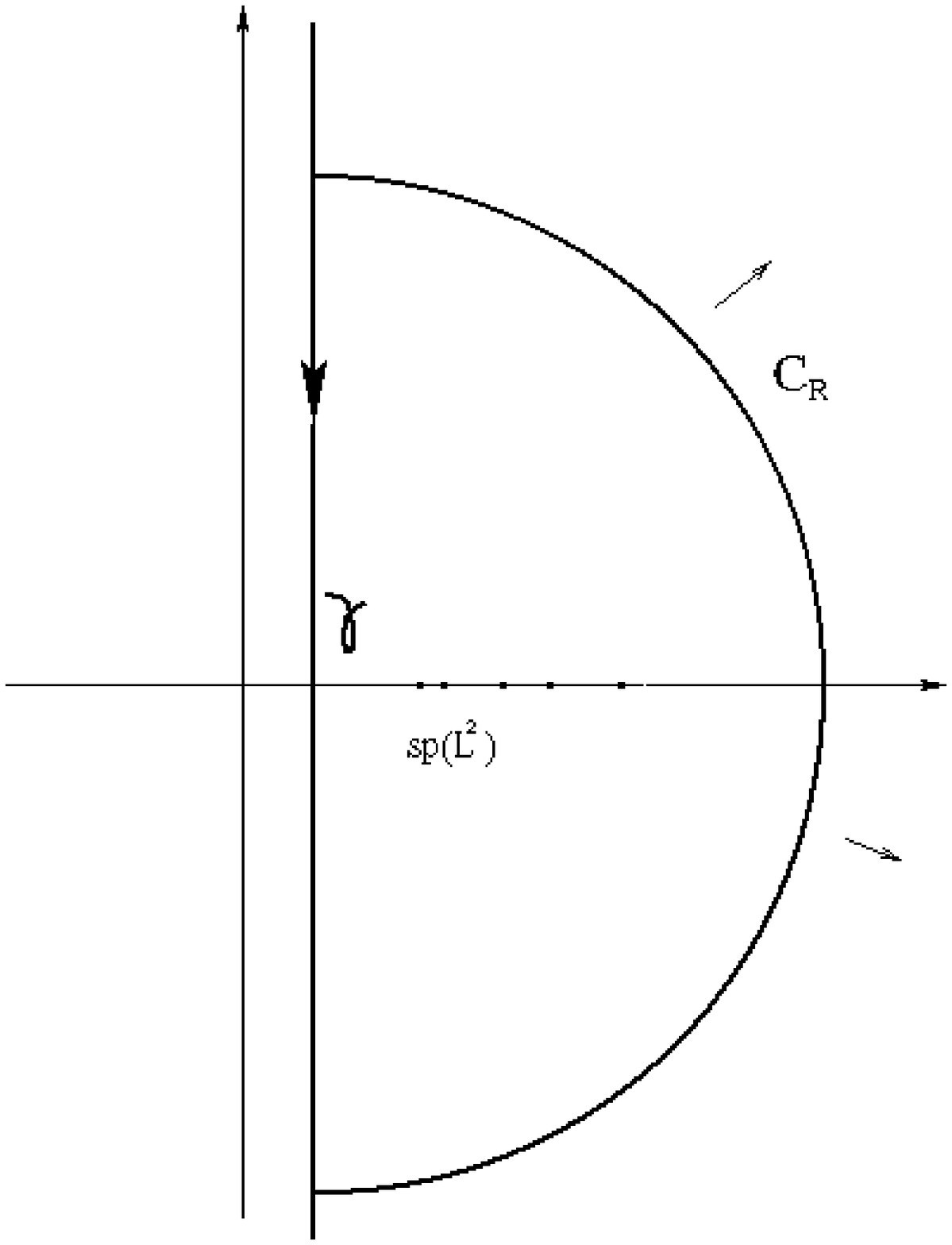,width=6cm, height=7cm}}
\caption{}
\end{figure}

Second, working on a single coordinate chart,
 $V_\alpha\bigcap
X_R$, we shall prove that $\nabla^2_L$ is polynomial in $\rho$.
Namely, for certain $G_0, G_1$ and $G_2,$ which are independent of
$\rho$, we show that
\begin{align}\label{section:details12}
\nabla^2_L=\rho^2G_2+\rho G_1+G_0.
\end{align}

Indeed, if we let $\nabla=d+\theta$, then:
\begin{align}\label{section:seppows}
\nabla_L^2&=\nabla^2+[\nabla,L]+L^2\\
\nonumber &= d\theta+\theta^2+[d+\theta,L]+L^2\\
\nonumber &=\underbrace{L^2}_{\rho^2G_2}+
  \underbrace{(d_\Xi L+d_xL+[\theta,L])}_{\rho G_1}
 +\underbrace{(d_x\theta+\theta^2+d\rho L)}_{G_0}.
\end{align}

Observe that  $G_0$ and $G_1$ are nilpotent of degree at most
$2n$, while $G_2$ is just $L^2/\rho^2$, hence invertible on $X_R$.

Third, by nilpotence of $G_0$ and $G_1$, the integrand
in (\ref{section:details0}) can be expanded in a terminating
geometric series:
\begin{align}\label{section:details1}
(-\nabla_L^2)^{-z}&={\frac 1{2\pi i}}\int_\gamma
{\scriptstyle \lambda^{-z}
  (\lambda +\rho^2G_2+\rho G_1 +G_0)^{-1}}d\lambda\\ \nonumber &
={\frac 1{2\pi i}}\int_\gamma
{\scriptstyle\lambda^{-z}(\lambda +\rho^2G_2)^{-1}}
\sum_{k=0}^{2n}(-1)^k{\scriptstyle \big[(\lambda +\rho^2G_2)^{-1}
(\rho G_1 +G_0)\big]^k}d\lambda.
\end{align}

Fourth, we separate the powers in  the $k$-th term of these series
and for certain sections $\Theta_l^k$ obtain:
\begin{align}\label{section:details4}
{\scriptstyle
(-1)^k(\lambda+\rho^2G_2)^{-1}
\big[(\lambda+\rho^2G_2)^{-1}
(\rho G_1+G_0)\big]^k=\sum_{l=0}^k\rho^{-2k+l}
\Theta^k_l(\lambda/\rho^2)d\rho d\Xi + err}.
\end{align}
The details of this computation are postponed until the end of the
proof. Here, the sections $\Theta^k_l$ will depend on the quantity
$\lambda/\rho^2$, but otherwise will not depend on $\rho$
explicitly. Also, $err$ represents the terms which may be ignored
because they are not multiples of the vertical volume form
$d\rho\,d\Xi$ and hence do not contribute to the integral over
$X_R$ in (\ref{section:details00}). Then, ignoring those error
terms, (\ref{section:details1}) becomes:
\begin{align}\label{section:details2}
(-\nabla^2_{L})^{-z}=\sum_{k,l}{\frac {\rho^{-2(z+k)+l}}{2 \pi i}}
\left[\int_{\gamma/\rho^2}\left({\frac\lambda{\rho^2}}\right)^{-z}
\Theta^k_l\left({\frac\lambda{\rho^2}}\right)
{\frac {d\lambda}{\rho^2}}\right]d\rho d\Xi.
\end{align}
Since the integrals can be computed  using the substitution
$\sigma=\lambda/\rho^2$
we can define $ \Theta^K(z,x,\Xi)$ as follows:
$$ \Theta^K(z,x,\Xi)=\sum_{l-2k=K}
      {\frac1 {2\pi i}}\int_{\gamma/\rho^2}\sigma^{-z}
      \Theta^k_ld\sigma.$$
Finally, observe that the above integrals do not really depend on
$\rho$ because the contour $\gamma/\rho^2$ can be replaced by  a
certain contour $\gamma'$ independent of $\rho$. This basically
finishes the proof, except we need to supply the following
details:
\begin{itemize}

\item[a)] The construction of $\Theta^k_l$ and the derivation
     of (\ref{section:details2}).
\item[b)] The exact choice of $\gamma'$.
\end{itemize}

For (a), we need to work out the details of
(\ref{section:details2}). Fix $\rho>0$. Consider some multiindex
$I_{k,l}=(\iota_1,\iota_2,\ldots,\iota_k)\in \{0,1\}^k$, in which
1 appears $l$ times and 0 appears $k-l$ times. Let $G'_0=G_0$ and
$G'_1=G_1\rho$. Then we define $G_{I_{k,l}}$ as
$$G_{I_{k,l}}=(\lambda+L^2)^{-1}\prod_{\mu=1}^k
\Big((\lambda+L^2)^{-1}G'_{\iota_\mu}\Big).$$ Then, recalling that
$\sigma=\lambda/\rho^2$ and $G_2=L^2/\rho^2$:
\begin{align}\label{section:weget1}
G_{I_{k,l}}&= (\lambda+L^2)^{-1}\prod_{\mu=1}^k
\Big((\lambda+L^2)^{-1}G'_{\iota_\mu}\Big)\\ \nonumber &=
\rho^{-2k+l} \big({\frac \lambda
{\rho^2}}+G_2\big)^{-1}\prod_{\mu=1}^k \Big(({\frac \lambda
{\rho^2}}+G_2)^{-1}G_{i_\mu}\Big)
\\ \nonumber
&= \rho^{-2k+l} \big(\sigma+G_2\big)^{-1}\prod_{\mu=1}^k
\Big((\sigma+G_2)^{-1}G_{i_\mu}\Big).
\end{align}

We define ${\tilde \Theta}^k_l$ as the sum of  those $G_{I_{k,l}}$
which are multiples of the vertical volume form $d\rho
d\Xi^1\ldots d\Xi^{n-1}=\rho^{n-1}d\rho d\Xi$. (All the others do
not contribute to the integral over $X_R$ in (\ref{theorem:1eq}) and
we ignore them).

Therefore we can pull both $\rho^{-2k+l}$ and $d\rho d\Xi$ out of
${\tilde \Theta}^k_l$ and define $\Theta^k_l$ through the
equation:
$${\tilde \Theta}^k_l=\rho^{-2k+l}\Theta^k_ld\rho d\Xi.$$
Now both (\ref{section:details4}) and (\ref{section:details2}) become
clear. The only possible issue is that as we change the variable
of integration from $\lambda$ to $\sigma$ in the Cauchy integral,
the contour of integration shifts:
\begin{align}\label{section:weget2}
\int_\gamma\scriptstyle\lambda^{-z}G_{I_{k,l}}d\lambda&=
\int_\gamma\scriptstyle
\lambda^{-z}(\lambda+L^2)^{-1}\prod_{\mu=1}^k
\Big((\lambda+L^2)^{-1}G'_{\iota_\mu}\Big)d\lambda\\ \nonumber
&=\int_{\gamma/\rho^2} \scriptstyle(\rho)^{-2(z+k)+l}\big({\frac
\lambda {\rho^2}}\big)^{-z} \big({\frac \lambda
{\rho^2}}+G_2\big)^{-1}\nonumber\\ \nonumber
&\qquad\qquad\,\,\scriptstyle\times\prod_{\mu=1}^k \Big(({\frac
\lambda {\rho^2}}+G_2)^{-1}G_{i_\mu}\Big) d\big({\frac \lambda
{\rho^2}}\big)\\
\nonumber &=\int_{\gamma'}
\scriptstyle(\rho)^{-2(z+k)+l}\sigma^{-z}
\big(\sigma+G_2\big)^{-1}\prod_{\mu=1}^k
\Big((\sigma+G_2)^{-1}G_{i_\mu}\Big) d\sigma.
\end{align}

Part (b)takes care of this issue. We choose $\gamma'$ to be the vertical
contour in ${\C}$, which is parameterized by $T-i\chi$,
($\chi\in{\R}$) for a suitable $T$.  We are about to show that
there exists $T$ such that
$$0<T<\inf\bigcup_{p\in X_R}sp(-\nabla_L^2(p)).$$
Here, the notation $\nabla_L^2(p)$ reminds us that
$\nabla_L^2$ is a section of $End(\pi^*E)$ which depends
on $p\in T^*M$, and $sp$ denotes the spectrum over each point $p$.

\begin{lemma} \label{section:Texists}
There exists $T>0$ as above.
In fact, there is an open subset $U$ of
$\{\lambda|Re(\lambda)>T\}$,
such that the pointwise spectrum $sp(-\nabla^2_L(p))$
is contained in $U$ for all $p\in X_R$.
\end{lemma}
{\bf Proof:} Indeed,
$\nabla^2_L$ equals $L^2$ plus the nilpotent term
$[\nabla,L]+\nabla^2$. So, $\lambda+\nabla_L^2$ is invertible
whenever $\lambda+L^2$ is invertible. This is apparent from
the geometric series (\ref{section:details1}).
 Thus, $sp(-\nabla^2_L)\subseteq sp(-L^2)$,
so it is enough to find $T$ such that
$$0<T<\inf\bigcup_{p\in X_R}sp(-L^2).$$
But by compactness of $S^*M$, there exists $T$ such that:
$$0<T<\inf\bigcup_{p\in S_R^*M}sp(-L^2).$$
Appealing to  homogeneity of $L^2$ in $\rho$, one sees that
$T$ satisfies the assertion of the lemma. This
finishes the proof of the lemma and the proposition.\qed

Elaborating on this proof, and {\it retaining the notation
therein}, we can obtain an estimate which shall be useful later:
\begin{lemma}\label{section:estim}
Suppose $\eta$ is compactly supported in a single chart
$U_\alpha$ of $M$.
Then there exist constants $K,K'$ such that for every
local section $\Theta^k_l$ of $\Lambda^*T^*M
\grtimes End(E)$ as in the proof of proposition
\ref{section:details}, and for all complex $z$:
\begin{align}
\left|\int_{S^*M\bigcap V_\alpha}
\pi^*(\eta)\operatorname{tr_s}
\int_{\gamma'}\sigma^{-z} \Theta^k_l(p)d\sigma\right|
\le Ke^{-K'Re(z)}
\end{align}
\end{lemma}
{\bf Proof:} By compactness of $S^*M$, there exists a closed loop
$\gamma''$ which simultaneously surrounds all the pointwise
spectra of $-L^2(p)$ for all $p\in S^*M$. Such a loop may be
chosen strictly to the right of the imaginary axis. This loop can
be used in place of the contour $\gamma'$ in the above integral
without affecting its value. The advantage is that $\gamma''$ is
compact. Then $|\sigma^{-z}|\le Ke^{-K'Re(z)}$ for all $\sigma\in
\gamma''$ and some constants $K, K'$, Also,
$\pi^*(\eta)\Theta^k_l$ is bounded on the compact set
$\gamma''\times S^*M\bigcap \pi^{-1}supp(\eta)$ by some $K''$
(with the appropriate choice of charts $U_\alpha$ and $V_\alpha$,
as assumed). Then, integrating out the variables $\sigma$, $x$,
and $\Xi$ over this set we see that:
\begin{align}\scriptstyle
\Big|\int_{S^*M\bigcap V_\alpha}
\pi^*(\eta)\operatorname{tr_s}
\int_{\gamma'}\sigma^{-z}&\scriptstyle \Theta^k_l(p)d\lambda\Big|
 \\ \nonumber
&\scriptstyle\le \int_{S^*M\bigcap V_\alpha} \int_{\gamma'}
\Big|\operatorname{tr_s}
\pi^*(\eta)\sigma^{-z}\Theta^k_l(p)\Big|d\lambda\\
\nonumber&\scriptstyle\le K'''K''Vol\Big((S^*M\bigcap
V_\alpha)\times\gamma''\Big)Ke^{-K'Re(z)},
\end{align} where $K'''$ comes from the supertrace.
\qed

\begin{corollary}\label{section:estcorollary}
For {\bf all} $\eta\in \Omega^*M$, there exist constants
$K$ and $K'$ such that:
$$I(z, \eta)<Ke^{-K'Re(z)}.$$\
\end{corollary}

{\bf Lemma \ref{theorem:poles}.}{\it
 Let $\eta\in\Omega^{\kappa}(T^*M)$.
Then $I(z,\eta)$ can only have a nonzero residue at the point
$(\kappa-2n)/2$. Further, if $n$ is even and the residue is
nonzero, then $(\kappa-n)/2$ must be an integer. If $n$ is odd,
$(\kappa-2n)/2$ must be a half-integer. In either case,
nonzero residues occur only for even $\kappa$.}\\
{\bf Proof:} Utilizing the proof of proposition
\ref{section:details} we reason out the case of even $n$, the odd
case being treated similarly.

Let $\eta$ be a $\kappa$-form. Recall that in in the said
proposition,
$$I(z,\eta)=\sum_{\alpha,k,l}\int_{X_R\bigcap V_\alpha}
\pi^*(\eta)\operatorname{tr_s}\int_\gamma
\sigma^{-z}\rho^{-2(z+k)+l}\Theta^k_l d\sigma.$$ We restrict our
attention to a single chart $V_\alpha$. We have seen that
$\Theta^k_l$ is expanded into the sum of terms $G_{I_{k,l}}$ using
(\ref{section:weget2}). The only such terms that could possibly
contribute to $I(z,\eta)$ are the ones of differential-form degree
$2n-\kappa$, because $\eta$ multiplied by them must produce a
$2n$-form on $T^*M$. Thus, we need to collect all the appropriate
(i.e., {\it contributing}) terms of the form:
\begin{align}\label{section:contrib}
\pi^*(\eta)\int_\gamma\lambda^{-z}G_{I_{k,l}}d\lambda.
\end{align}

The remainder of the proof is but an exercise in counting the
differential form degrees. They are products of $l$ copies of
$G_1$ (which is locally a matrix of 1-forms), and $k-l$ copies of
$G_0$, which is $\nabla^2+d_\rho L$. Thus, any contributing term
requires $k-l\ge1$, because at least one copy of $G_0$ is needed
to supply the differential $d\rho$ for the $2n$-form. Each
additional copy of $G_0$ can only supply the curvature $\nabla^2$,
which is a matrix of 2-forms. Therefore, each additional copy of
$G_0$ may be replaced by two copies of $G_1$ without changing the
total degree of (\ref{section:contrib}). So, all the contributing
terms satisfy:
\begin{align}\label{section:contrib2}
 deg(\nabla^2)(k-l-1)+deg(d_\rho L)+deg(G_1)l=2n-\kappa,
\end{align}
which means that the quantity $2k-l$ is the same for all of them.
But  one can see from (\ref{section:details2}) that
 the location of the residue which arises from
each contributing term of (\ref{section:weget2})
depends only on that quantity.
It is apparent from (\ref{section:weget2}) that
 $I(z,\eta)$ can have at most one nonzero residue, whose
location must be ${\frac 12}(l-2k)={\frac {\kappa-2n}2}$.
This location does not depend on the topological information
about $M$ (other than its dimension). Neither does it depend
on the vector bundle $E$, on the curvature, etc.

If $\kappa$ is even,
then by (\ref{section:contrib2}), $l$ must be odd.
Similarly, if $\kappa$ odd, then $l$ must be even.
Now, suppose we equip $\Lambda^*T^*M\grtimes End(\pi^*E)$
with the ${\Z}_2$-grading inherited from $End(\pi^*E)$,
rather than the total one. The supertrace vanishes
on the sections of $\Lambda^*T^*M\grtimes End(\pi^*E)$
which are odd in that inherited grading. We
say that such sections are {\it of an odd profile}.
The term {\it even profile} is defined similarly.

Hence, if $l$ is even, we get a term
$\int_\gamma\lambda^{-z}G_{I_{k,l}}d\lambda$ which
involves some number of the even-profile sections
$(\sigma+L^2)^{-1}$ and $\nabla^2$. Also, it involves
$d_\rho L$ and an even number of copies of $G_1$.
Such a term will be of odd profile and its
supertrace is zero. Thus,  only if $\kappa$ is even can
one hope to get a non-zero residue. Combining
(\ref{section:contrib2}) with the fact that the
residue is located at $(l-2k)/2$, we get
$l-2k=-2n+\kappa$.
\qed

\section{Proof of Theorem \ref{theorem:2}.}\label{section:pt2}
{\bf Theorem \ref{theorem:2}}\begin{it}
Under the hypotheses outlined in section
\ref{section:Form}, for any  $\eta\in\Omega^*(M)$
\begin{multline*}
\int_{T^*M}\operatorname{tr_s}\pi^*(\eta)\exp{\nabla_L^{2}} =
\lim_{R\to 0}\sum_{z\in{\C}}Res|_z\Gamma(z)
\int_{X_R}\operatorname{tr_s}\pi^*(\eta)({-\nabla_L^2})^{-z}.
\end{multline*}
Further, all but finitely many residues on the
right-hand side vanish as $R\to 0$.
\end{it}

In fact, we shall prove that for any $R>0$:

$$
\int_{X_R}\pi^*(\eta)\operatorname{tr_s}\exp\nabla^{2}_{L} =
\sum_{z\in{\C}}Res|_z\Big[\Gamma(z)
\int_{X_R}\pi^*(\eta)\operatorname{tr_s}(-\nabla_{L}^2)^{-z}\Big].
$$

Taking limits of both sides as $R\to 0$, we get  theorem
\ref{theorem:2}.
The  outline of our proof is the following.

\begin{itemize}
\item[1)]
We introduce a parameter $t\ge0$ and express $\exp\nabla_{tL}^{2}$
using a Cauchy integral.
Thus, for a suitable vertical contour $\gamma$ in ${\C}$,
\begin{align}\label{section:bad}
\exp\nabla_{tL}^{2}=
\int_\gamma e^{-\lambda}(\lambda+\nabla_{tL}^2)^{-1}
d\lambda.
\end{align}
Strictly speaking,
(\ref{section:bad}) is incorrect, because the integral
over $\gamma$ does not converge. Still, the formula
holds in a weak sense. That is,
for all $R>0$ and $\eta\in\Omega^*M,$
\begin{align}\label{section:ol1}
\int_{X_R}\pi^*(\eta)\exp\nabla_{tL}^{2}=
\int_{X_R}\int_\gamma \pi^*(\eta)e^{-\lambda}
(\lambda+\nabla_{tL}^2)^{-1}d\lambda,
\end{align}
where the integral over $\gamma$ converges.
This is proven using the geometric-series expansion of
$e^{-\lambda}(\lambda+\nabla_{tL}^2)^{-1}$ very similar to
that in (\ref{section:details00}). We show that the above integral
$\int_\gamma$ converges at least for those terms of the
expansion which do contribute to
(\ref{section:ol1}).  See lemma \ref{theorem:poles} for
discussion of contributing and non-contributing terms.
Thus, (\ref{section:ol1})
is an equality of {\it currents} on
$\Omega^*M$ induced by $\operatorname{tr_s}\exp\nabla^2_{tL}$
and by $\int_\gamma e^{-\lambda}
(\lambda+\nabla_{tL}^2)^{-1}d\lambda$ via
 Poincare duality, as discussed in section \ref{section:Form}.
\item[2)]
Assuming $Re(z)\gg1$, we show that in the same weak sense,
\begin{align}
\int^\infty_0t^{z-1}\exp\nabla_{tL}^{2}dt=
\Gamma(z)\left(-\nabla^{2}_{L}\right)^{-z},
\end{align}
which means that:
\begin{align}\label{section:ol2}
\int_{X_R}\pi^*(\eta)\int^\infty_0t^{z-1}\exp\nabla_{tL}^{2}dt=
\Gamma(z)\int_{X_R}\pi^*(\eta)\left(-\nabla^{2}_{L}\right)^{-z}.
\end{align}
This is an  application of the so-called {\it Mellin transform}
which is is given by $f\mapsto \int_0^\infty f(t)t^{z-1}dt$.
The inverse transform is given by
 $F(z)\mapsto{\frac 1 {2\pi i}}\int_Ct^{-z}F(z)dz$,
where $C$ is a suitable vertical contour ${\C}$.
See \cite{Arfken} for details.
Roughly, the computation
for (\ref{section:ol2}) is the following:
\begin{align}\label{section:ol3}
\int_0^\infty\int_{X_R}&\pi^*(\eta)\operatorname{tr_s}
\exp\nabla_{tL}^{2}t^{z-1}dt\\&=
\int_{X_R}\pi^*(\eta)\operatorname{tr_s} \int_0^\infty t^{z-1}
\int_\gamma e^{-\lambda}(\lambda+\nabla_{tL}^2)^{-1}d\lambda dt
\nonumber\\
&=
\int_{X_R}\pi^*(\eta)\operatorname{tr_s}
\int_\gamma \int_0^\infty
e^{-\lambda}(\lambda+\nabla_{tL}^2)^{-1}t^{z-1}dtd\lambda \nonumber\\
&=
\int_{X_R}\pi^*(\eta)\operatorname{tr_s}
\int_\gamma
\Gamma(z)\lambda^{-z}(\lambda+\nabla_{L}^2)^{-1}d\lambda \nonumber\\
&=\Gamma(z)
\int_{X_R}\pi^*(\eta)\operatorname{tr_s}
(-\nabla_{L}^2)^{-z}.\nonumber
\end{align}
Interchanging the integrals $\int_0^\infty\ldots dt$ and
$\int_{X_R}$ is easy, because $L^2$ is
negative definite and invertible, so $\exp\nabla^2_{tL}$
is rapidly decreasing and absolutely integrable on
${\R}\times X_R$.
We will need to prove that we can
interchange $\int_0^\infty\ldots dt$ and $\int_\gamma\ldots d\lambda$,
at least for $Re(z)\gg1$ and for the relevant terms of the
geometric series expansion of $e^{-\lambda}(\lambda+\nabla_{tL}^2)^{-1}$.
Also, notice that the exponential is $e^{-\lambda}$, not
$e^{-t\lambda}$, as one might expect from the well-known identity:
$$\lambda^{-z}\Gamma(z)=\int^\infty_0e^{-t\lambda}t^{z-1}dt.$$
The algebra behind this will be explained.
Also, observe that by theorem \ref{theorem:1} the
integral $\int_{X_R}\ldots$ in the right-hand
side has a meromorphic extension to ${\C}$ with at
most simple poles.
\item[3)]
Next, we simply restate (\ref{section:ol3})
in terms of the inverse Mellin transform.
For a certain vertical $C\subset {\C}$,
\begin{align}\label{section:ol4}
\int_{X_R}\!\!\!\!\!
\operatorname{tr_s}\pi^*(\eta)\exp\nabla^{2}_{tL} =
{\frac 1{2\pi i}}\int_C\!\!t^{-z}\Gamma(z)
\int_{X_R}\!\!\!\!\operatorname{tr_s}
\pi^*(\eta)(-\nabla_{tL}^2)^{-z}dz.\!\!\!\!\!
\end{align}
The integral $\int_{X_R}\!\!\!\!\operatorname{tr_s}
\pi^*(\eta)(\!-\!\nabla_{tL}^2)^{-z}\!dz\!$
will be abbreviated by $I_t(z,\eta)$.
Since the meromorphic extension of $I_t(z,\eta)$
is defined for all $z$ except at a certain discrete set,
we can choose $C$ to pass through
the domain where $I_t(z,\eta)$ is defined.
We will need to check the convergence of the integral
$\int_C\ldots dz$. Observe that $t$ has reappeared
in the subscript $tL$ on the right-hand side.
This will require clarification.
\item[4)]
After taking the meromorphic extension of $I_t(z,\eta)$
in the right-hand side of (\ref{section:ol4}),
 we ``collect'' the residues
by translating the contour $C$ to the left.
We will need to prove that as the contour translates,
$$\int_C\!\!\Gamma(z)\!
\int_{X_R}\!\!\!\!\operatorname{tr_s}\pi^*\!(\eta)(-\nabla_{tL}^2)^{-z}\!dz
\!\!\longrightarrow\!\!\sum_{z\in{\C}}Res|_z \!\Gamma(z)
\int_{X_R}|!\!\!\!\operatorname{tr_s}\pi^*(\eta)(-\nabla_{tL}^2)^{-z}.$$
This is an application of the residue theorem (Fig. 2).
\end{itemize}

\begin{figure}
\centerline {\psfig {figure=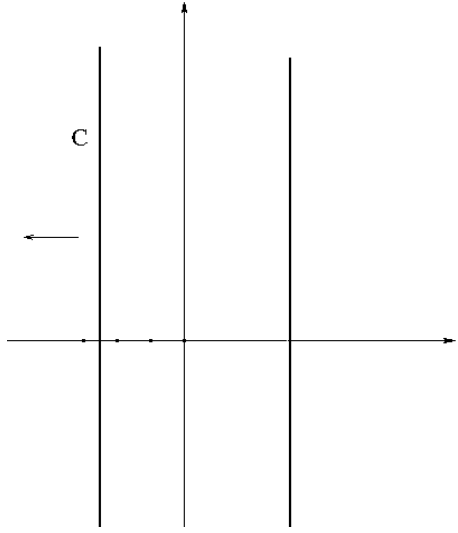,width=6cm, height=7cm}}
\caption{Collecting the residues.}
\end{figure}

We proceed with the proof. In order to make the equation
(\ref{section:ol1}) in step (1) less bulky, we  introduce
the following notation.

\begin{defn}\label{section:weak}
Let $\mu, \nu$ be smooth forms on $T^*M$. We say that they are
{\it equal in the weak sense} and we write $\mu=_w\nu$ if
for all $\eta\in \Omega^*(M)$
$$\int_{X_R}\pi^*(\eta)\mu
=\int_{X_R}\pi^*(\eta)\nu.$$
In other words, this means that $\mu$ and $\nu$ are equal
as currents on $\Omega^*(M)$.
\end{defn}
With this kind of equivalence relation, it is possible to
write $\exp(\nabla^2_{tL})$ as a Cauchy integral, similarly to
(\ref{section:details0}).
But first,  we introduce some more notation.
We expand $\nabla^2_{tL}$, similarly to (\ref{section:details1}):
$$\nabla^2_{tL}=\nabla^2+t[\nabla,L]+t^2L^2.$$
Now, let $\Theta_t=\nabla^2+t[\nabla,L]$,
and let $\gamma_\tau$  denote a vertical
line $\{Re(\lambda)=\tau\}$, oriented downward. (E.g.,
$\gamma_0$ is just the imaginary axis).
Since locally, $\Theta_t$ is a matrix of differential forms
of positive degree, and $T^*M$ is $2n$-dimensional, we
have the following geometric series:
\begin{align}\label{section:gseries}
(\lambda+\nabla_{tL}^2)^{-1}=
(\lambda+t^2L^2)^{-1}\sum^{2n}_{k=0}
(-1)^k\big[(\lambda+t^2L^2)^{-1}\Theta_t\big]^k.
\end{align}

\begin{lemma}\label{section:tech3}
Fix $p\in T^*M$.
Choose $\epsilon\ge 0$ such that
$\gamma_\epsilon$ lies to the left of the {\it pointwise}
spectrum $sp(-t^2L^2(p))$ of $-t^2L^2(p)$.
Then  the following equality holds and the integral
on the right-hand side converges:
\begin{align}\label{section:tech3eq}
\exp(&\nabla_{tL}^2(p))\!=_w\\ \nonumber &{\frac 1 {2\pi
i}}\!\!\int_{\gamma_\epsilon}
\!\!\!\!{e^{-\lambda}\!\sum_{k=1}^{2n}\!
(\lambda+t^2L^2(p))^{\!\!^{-1}}\!\!
(-1)^k\big((\lambda+t^2L^{2}(p))^{^{-1}}\! \Theta_t\big)^k }
\!d\lambda.
\end{align}
\end{lemma}
In view of (\ref{section:gseries}), this comes close to
\begin{align}
\exp(\nabla_{tL}^2)={\frac 1 {2\pi i}}\int_{\gamma_\epsilon}
 e^{-\lambda} (\lambda+\nabla^2_{tL})^{-1}
d\lambda.
\end{align}

{\bf Proof:}
The contributing terms (see the proof of lemma
\ref{theorem:poles}) have enough negative powers of $\lambda$
to assure convergence.
Note that we start the series at $k=1$ which explains the weak
equality: the integral of the term with $k=0$
diverges. Fortunately,
that term does not involve vertical
differentials $d\Xi^i$ or $d\rho$. Thus, by the proof of lemma
\ref{theorem:poles}, it does not contribute to any current.
\qed

In step (2), we need to prove some version of (\ref{section:ol2})
 (it is not true literally).
To do this, we shall:
\begin{itemize}
\item [a)]Apply the geometric-series expansion (\ref{section:gseries}) to
$$e^{-\lambda}(\lambda+\nabla_{tL}^2)^{-1}=
t^2L^2+t[L,\nabla]+\nabla^2.$$
\item[b)] On the $k$-th term of that expansion, perform a secondary
         expansion into terms which will be denoted by $\Phi^k_l$.
         This will separate the powers of $t$:
\begin{align}\label{section:mellin1}
(-1)^k(\lambda+t^2L^2)^{-1}\left[(\lambda+t^2L^2)^{-1}\Theta_t\right]^k=
\sum_{l=0}^kt^{l-2(k+1)}\Phi^k_l.
\end{align}
\item [c)]For each term $\Phi^k_l$, take the Mellin transform
          of the integral $
       \int_{\gamma}e^{-\lambda}t^{l-2(k+1)}\Phi^k_ld\gamma$.
          Homogeneity in $t$ makes this task possible.
          This will yield a formula
          much like (\ref{section:ol2}),  for each
          of a collection of separate terms.
          At a certain later point, we shall reassemble the Mellin
          transforms of these terms into the quantity
          $\Gamma(z)(-\nabla^2_L)^{-z}$.
\end{itemize}
We accomplish a), b) and c) in
the following proposition.  Also, our earlier warning about
interchanging the integrals $\int^\infty_0\ldots dt$
and $\int_\gamma\ldots d\lambda$ in (\ref{section:ol2})
receives due attention here.

\begin{proposition}\label{section:mellin}
There exist smooth sections $\Phi^k_l$
of\\ $\Lambda^*T^*M\grtimes End(\pi^*E)$ such that the
following expansion holds for each $k>0$:
\begin{align}
(-1)^k(\lambda+t^2L^2)^{-1}\left[(\lambda+t^2L^2)^{-1}
\Theta_t\right]^k=
\sum_{l=0}^kt^{l-2(k+1)}\Phi^k_l.
\end{align}
These sections  depend on the quantity
${\frac \lambda {t^2}}$, but they do not depend on $t$
in any other way.
Further,  there exists $\tau>0$ such that for $Re(z)\gg0$,
\begin{align}\label{section:mellin2}
\int_0^\infty &t^{2z-1}\operatorname{tr_s}\exp\nabla_{tL}^2dt\\
&=_w\int_0^\infty t^{2z-1}\operatorname{tr_s}
  \Big[\int_{\gamma_{_0}}\sum_{k=1}^{2n}
    e^{-\lambda}
    (-1)^k(\lambda+t^2L^2)^{-1}\nonumber\\
         &\qquad\qquad\quad\qquad\qquad\qquad\times\big[(\lambda+t^2L^2)^{-1}\Theta_t\big]^k
    d\lambda \Big]dt\nonumber\\
&={\frac 1 2} \int_{\gamma_{\tau}}\operatorname{tr_s}
\sum_{k=1}^{2n}\sum_{l=0}^{k}
    \Gamma(z+l/2-k)\Big({\frac \lambda{t^2}}\Big)^{-(z+l/2-k)}
        \Phi_l^k{\frac {d\lambda}{t^2}}.\nonumber
\end{align}
\end{proposition}
Here, $\gamma_{_0}$ is the
imaginary axis. Also, note that we are using $2z$ instead of $z$
in the Mellin transform.

{\bf Proof:}
This proof is very similar to that of proposition
\ref{section:details}.
Each term $\big[(\lambda+t^2L^2)^{-1}\Theta_t\big]^k$ is a
non-commutative polynomial in the quantities $\nabla^2$ and $t[L,\nabla]$,
which are globally well-defined sections of the vector bundle
$\Lambda^*T^*M\grtimes End(\pi^*E)$.
Therefore, each term can be further further expanded as follows:
$$(\lambda+t^2L^2)^{-1}
(-1)^k\big[(\lambda+t^2L^2)^{-1}\Theta_t\big]^k=\sum_l{\tilde \Phi^k_l},$$
where the quantity ${\tilde \Phi^k_l}$ is the sum of all
the monomials which are products of $l$ copies of $t[L,\nabla]$,
 $k-l$ copies of $\nabla^2$, and $k+1$ copies of $(\lambda+t^2L^2)^{-1}$.
This is similar to the construction of $\Theta^k_l$ in
the proof of proposition \ref{section:details}.
But $(\lambda+t^2L^2)^{-1}=t^{-2}({\frac\lambda{t^2}}+L^2)^{-1}$,
so, we can pull $t^{-2(k+1)+l}$ out of
${\tilde\Phi^k_l}$ to obtain $\Phi^k_l$:
\begin{align}\label{section:defphi}
\tilde\Phi^k_l=t^{l-2(k+1)}\Phi^k_l.\end{align}

So (\ref{section:mellin1}) holds and we can prove (\ref{section:mellin2}.)
Just as in the lemma \ref{section:tech3}, the weak equality
is there because we start the geometric series at $k=1$.
In what follows, by $sp(t^2L^2)$ we always mean the
pointwise spectrum over a point $p\in T^*M$. Define $\tau$ by
$$\tau={\frac 12}\inf\bigcup_{p \in X_R}sp(-L^2(p)).$$
Such $\tau$ exists by lemma \ref{section:Texists}.
Then for each $t>0$, the pointwise spectrum of $t^2L^2$ is to the
right of $\gamma_{\tau t^2}$.  Fixing one such $t$ for the moment,
we see that the vertical $\gamma_{\tau t^2}$
is a suitable contour for the Cauchy integral
expression of $\exp\nabla^2_{tL}$
(by lemma \ref{section:tech3})
and we can compute, for the $k$-th term:
\begin{align}\label{section:st21}
\int_{\gamma_{ {\tau t^2} }}e^{-\lambda}(-1)^k
    (\lambda+t^2L^2)^{-1}
         \big[(\lambda&+t^2L^2)^{-1}\Theta_t\big]^k
    d\lambda \\\nonumber
&=\int_{\gamma_{\tau }}
 e^{-t^2\sigma}
    \sum_lt^{l-2k}\Phi_l^kd\sigma,
\end{align}
where $\sigma={\frac \lambda{t^2}}$. Thus,
\begin{align}\label{section:st22}
\int_0^\infty t^{2z-1}
  \Big[\int_{\gamma_{ {\tau t^2} }}&e^{-\lambda}
    (-1)^k(\lambda+t^2L^2)^{-1}
         \big[(\lambda+t^2L^2)^{-1}\Theta_t\big]^k
    d\lambda \Big]dt\\
&=\int_0^\infty t^{2z-1}
  \Big[\int_{\gamma_{\tau }}
 e^{-t^2\sigma}
    \sum_lt^{l-2k}\Phi_l^kd\sigma\Big]dt.\nonumber
\end{align}

Since $k>0$, $\Phi^k_l$ involves at least 2
factors of $(\sigma+L^2)^{-1}$. Therefore it is
absolutely integrable with respect to $\sigma$, while
$t^{2(z-k)+l-1}e^{-t^2\sigma}$ is absolutely integrable in $t$
for nonnegative $Re(z)$.
By Fubini's theorem, we may interchange the integrals
and finish the computation:

\begin{align}\label{section:st23}
\int_0^\infty t^{2z-1} &
  \Big[\int_{\gamma_{ {\tau t^2} }}e^{-\lambda}
    (-1)^k(\lambda+t^2L^2)^{-1}
         \left[(\lambda+t^2L^2)^{-1}\Theta_t\right]^k
    d\lambda \Big]dt\\\nonumber
&= \int_0^\infty t^{2z-1}
  \Big[\int_{\gamma_{\tau }}
 e^{-t^2\sigma}
    \sum_lt^{l-2k}\Phi_l^kd\sigma\Big]dt\\\nonumber
&=
 \int_{\gamma_{\tau  }}
  \Big[
\int_0^\infty t^{2z-1}
 e^{-t^2\sigma}
    \sum_lt^{l-2k}\Phi_l^kdt\Big]d\sigma\\\nonumber
&= {\frac 1 2} \int_{\gamma_{ \tau }}
  \Big[\int_0^\infty
 e^{-t^2\sigma}
    \sum_lt^{2z+l-2k-2}\Phi_l^kdt^2\Big]d\sigma\\\nonumber
&={\frac 1 2} \int_{\gamma_\tau}
    \sum_l\Gamma(z+l/2-k)\sigma^{-(z+l/2-k)}\Phi_l^kd\sigma.
\,\,\,\square\nonumber
\end{align}
This proves that
\begin{align}\label{section:ol3bisbis}
\int_{X_R}\!\!\pi^*(\eta)&\operatorname{tr_s} \int_0^\infty
t^{2z-1}\!\!\!\!
  \int_{\gamma_{ {\tau t^2} }}\! \!\!\!\!\!{e^{-\lambda}
    (-1)^k(\lambda+t^2L^2)^{-1}
         \big[(\lambda+t^2L^2)^{-1}\Theta_t\big]^k
}
    d\lambda dt\nonumber\\
&=\int_{X_R}\pi^*(\eta)\operatorname{tr_s}
    {\frac 1 2} \int_{\gamma_\tau}
    \sum_l\Gamma(z+l/2-k)\sigma^{-(z+l/2-k)}\Phi_l^kd\sigma,
\end{align}
and the integrals $\int_{X_R}$ and $\int_0^\infty\ldots\,dt$
in the left-hand side can be interchanged, as  remarked in our discussion
after (\ref{section:ol3}).

This equation is as close as we get to  (\ref{section:ol3}).
Our next step is to apply the inverse Mellin transform
to the right-hand side. Our estimate from lemma \ref{section:estim}
comes in here. The inverse Mellin transform
of (\ref{section:ol3bisbis}) is
\begin{align}
{\frac 1{2\pi
i}}\int_C&t^{-2z}\int_{X_R}\pi^*(\eta)\times\\\nonumber
 &\operatorname{tr_s}
  {\frac 1 2} \int_{\gamma_\tau}
  \sum_l\Gamma(z+l/2-k)\sigma^{-(z+l/2-k)}\Phi_l^kd\sigma\,dz.
\end{align}

By theorem \ref{theorem:1}, the vertical line $C$ may be chosen
very far to the right so the residues of $\Gamma(z)I_t(z,\eta)$
are nowhere near.
Convergence of the integral over $C$ is assured by
the estimate very similar to lemma \ref{section:estim}
and by the fact that $\Gamma(z)$ is rapidly decreasing on the
vertical lines.

\begin{align}\label{section:ol4bis}
\int_{X_R}\pi^*(\eta)&\operatorname{tr_s} \exp(\nabla_{tL})^2
\\\nonumber
 &= \sum_{k,l}
{\frac 1 {4\pi i}}\int_Ct^{-2z}\int_{X_R}\pi^*(\eta)  \\\nonumber
         &\qquad\times\Gamma(z+l/2-k)\operatorname{tr_s}
            \int_{\gamma_\tau}
              \sigma^{-(z+l/2-k)}\Phi_l^kd\sigma\, d(2z)\\\nonumber
&={\frac 1 {2\pi i}} \sum_{k,l}\int_C t^{-2z}
              \Gamma(z+l/2-k)\\\nonumber
             &\qquad\times
             \int_{X_R}\pi^*(\eta)\operatorname{tr_s}
                   \int_{\gamma_\tau}
              \sigma^{-(z+l/2-k)}\Phi_l^kd\sigma \,dz.
\end{align}
Introducing the variables $s=z+l/2-k$ and the verticals $C_{l,k}=C+l/2-k$,
we may rewrite (\ref{section:ol4bis}) as:
\begin{align}
{\frac 1 {2\pi i}} \sum_{k,l}\int_{C_{l,k}} t^{-2(s-l/2+k)}\Gamma(s)
                \int_{\gamma_\tau} \sigma^{-s}\Phi_l^kd\sigma ds.
\end{align}
In view of the next lemma, we may replace all the contours $C_{k,l}$
with $C$.

\begin{lemma} \label{section:moveC}
Fix some $p\in T^*M$.  Let $C$ and $C'$ be two vertical lines
in ${\C}$  with the same orientation.
If the expression $\Gamma(s)\int_\gamma \sigma^{-s}\Phi_l^kd\sigma ds$
has no singularities between them, then:
\begin{align}
\int_{C} t^{-2(s-l/2+k)}
        &\Gamma(s)\int_{\gamma_\tau}
         \sigma^{-s}\Phi_l^k(p)d\sigma ds=\\\nonumber
      & \int_{C'} t^{-2(s-l/2+k)}
                \Gamma(s)\int_{\gamma_\tau}
         \sigma^{-s}\Phi_l^k(p)d\sigma ds.
\end{align}
\end{lemma}
{\bf Proof:}
First, we  join $C$ and $C'$ by horizontal line segments
$ab$ and $cd$, located below and above the real axis,  as in Fig. 3.
Then:
\begin{align}
\left(\int_a^b\!\!\!+\int_b^c\!\!\!+\int_c^d
\!\!\!+\int_d^a\right)
   t^{-2(s-l/2+k)}\left[\int_{\gamma_\tau}
\!\!\!\Gamma(s)  \sigma^{-s}\Phi_l^kd\sigma\right] \ ds=0.
\end{align}

Because  on the vertical lines
$\Gamma(s)$ is rapidly decreasing
and the rest of the integrand is bounded in $s$,
the integrals over $ab$ and $cd$ vanish as those line segments
move away from the real axis. The result follows. \qed

\begin{figure}
\centerline {\psfig {figure=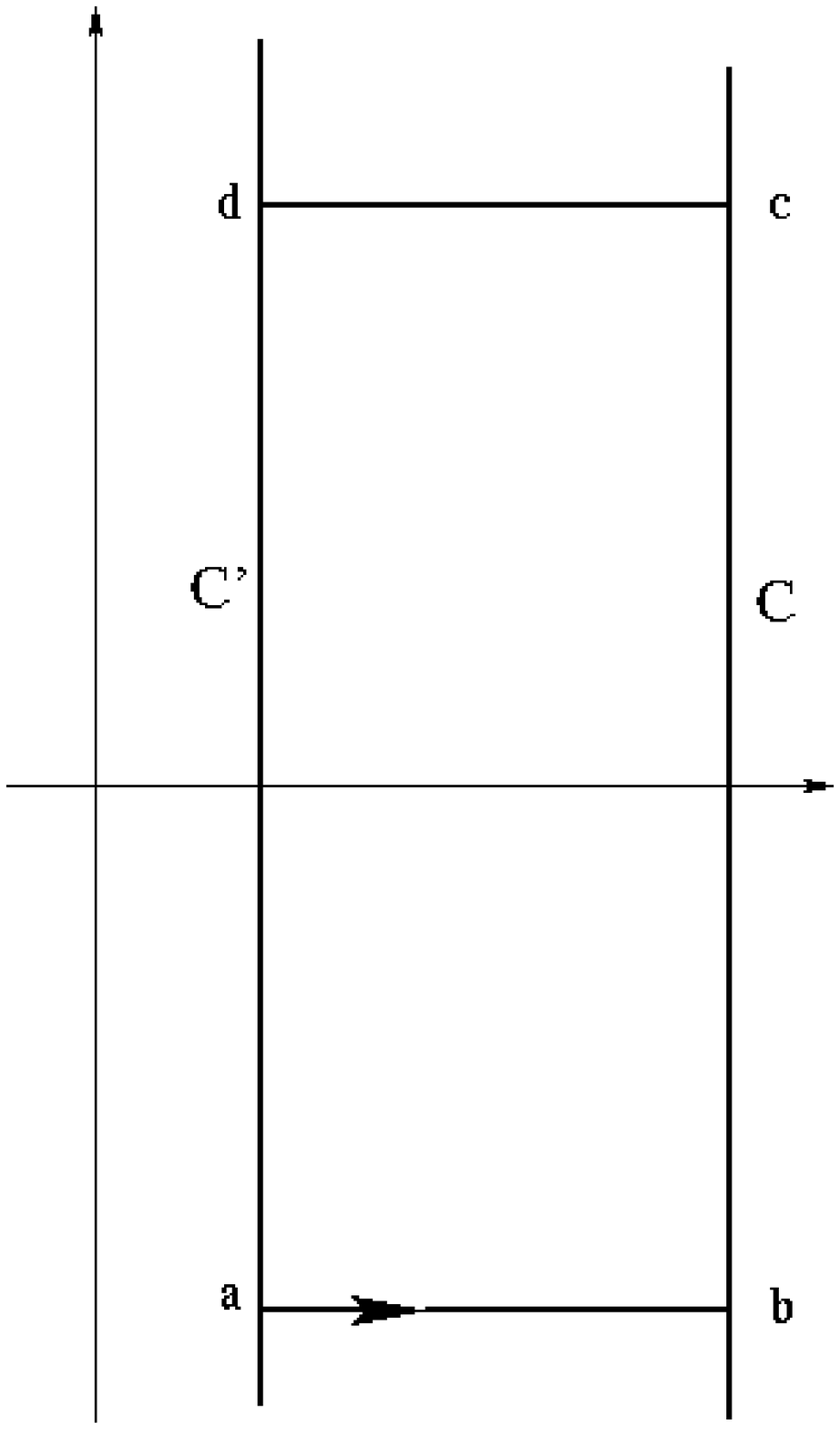, width=6cm, height=6cm}}
\caption{Proof of lemma \ref{section:moveC}.}
\end{figure}

This lemma allows us to continue the computation,
using $C$ instead of $C_{k,l}$, provided that $C$ were originally
chosen far enough to the right. We are about to reassemble
the individual terms $\int_{\gamma_\tau}\sigma^{-s}\Phi^k_ld\lambda$ into
the quantity $(-\nabla_L^{2})^{-z}$, as promised earlier.
\begin{align}
{\frac 1 {2\pi i}}\sum_{k,l}&\int_{C}
\int_{X_R}\pi^*(\eta)\operatorname{tr_s}
t^{-2(s-l/2+k)}\int_{\gamma_\tau}
              \Gamma(s)\sigma^{-s}\Phi_l^kd\sigma ds\\\nonumber
&={\frac 1 {2\pi i}} \sum_{k,l}\int_{C}
\int_{X_R}\pi^*(\eta)\operatorname{tr_s}
 t^{-2s}\int_{\gamma_\tau}
              \Gamma(s)\sigma^{-s}t^{l-2k}\Phi_l^kd\sigma ds.
\label{section:summingl}
\end{align}
Recall that $\Phi^k_l$ depends on the quantity
$\sigma={\frac \lambda {t^2}}$. Also, recall (\ref{section:defphi}):
\begin{align}
\tilde\Phi^k_l=t^{l-2(k+1)}\Phi^k_l.\end{align}
So, by (\ref{section:st21}), summing
the right-hand side of (\ref{section:summingl}) over $l$ we
obtain:

\begin{align}
{\frac 1 {2\pi i}}&\sum_{l}\!\int_{C}\!\!
\int_{X_R}\!\!\!\!\!\pi^*(\eta)\operatorname{tr_s}
t^{-2(s-l/2+k)}\int_{\gamma_\tau}\!\!\!\!\!
              \Gamma(s)\sigma^{-s}\!\Phi_l^kd\sigma ds\\\nonumber
&={\frac 1 {2\pi i}} \sum_{l}\int_{C}
\int_{X_R}\pi^*(\eta)\operatorname{tr_s}
 \int_{\gamma_{\tau t^2}}
          \Gamma(s)\lambda^{-s}
             {\tilde \Phi_l^k}d\lambda ds\\
&=\!\!{\frac 1 {2\pi i}} \int_{C}\!
\int_{X_R}\!\!\!\!\!\!\!\pi^*\!(\eta)\nonumber\!\times\!\!\operatorname{tr_s}\!\!
   \int_{\gamma_{ {\tau t^2} }}\!\!\!\!\!\!\!\lambda^{-z}(-1)^k
    (\lambda\!+t^2L^2)^{\!\!^{-1}}\\\nonumber&
         \!\qquad\qquad\qquad\qquad\qquad\qquad\times\big[(\lambda\!+t^2L^2)^{\!\!^{-1}}\!\Theta_t\big]^k\!
    d\lambda.
\end{align}
Finally, by summing this over all $k$, we recover the quantity
\begin{align}\label{section:summingk}
{\frac 1{2\pi i}}\int_C\int_{X_R}\pi^*(\eta)\operatorname{tr_s}
\Gamma(s)(-\nabla_{tL})^{-2s}ds.
\end{align}
We now invoke theorem \ref{theorem:1} and lift our standing
assumption that $Re(z)\gg0$. Hence, the
integral $\int_{X_R}\ldots$ in (\ref{section:summingk})
may be replaced by its meromorphic extension.
We can now finish the
proof, by moving the vertical $C$ to the left
and ``picking up'' all the residues. The procedure is explained
in the following  lemma.

\begin{lemma} For any $r>0$ such that
$(C-r)$ does not intersect the real axis at any of the residues,
the following holds:
\begin{align}\label{section:moveC1}
\int_{C}\Gamma(s)\int_{X_R}\pi^*(\eta)&(-\nabla_{tL}^2)^{-s}ds=\\\nonumber
&\int_{C-r}\Gamma(s)\int_{X_R}\pi^*(\eta)(-\nabla_{tL}^2)^{-s}ds\,
+ \\\nonumber &\,\,\, \sum_{Re(s)>-r} Res|_s
\Gamma(s)\int_{X_R}\pi^*(\eta)(-\nabla_{tL}^2)^{-s}.
\end{align}
Further,  substituting $r_m={\frac {2m+1}2}$ instead of $r$ in the
above expression,
the integral on the right-hand side tends to
zero as $m\to \infty$, so that
\begin{align}\label{section:moveC2}
\int_{C}\Gamma(s)\int_{X_R}\pi^*(\eta)&(-\nabla_{tL}^2)^{-s}ds=
\\\nonumber &\sum_{s\in {\C}} Res|_s
\Gamma(s)\int_{X_R}\pi^*(\eta) (-\nabla_{tL}^2)^{-s}.
\end{align}
\end{lemma}

{\bf Proof:}
The first equation follows by the argument similar to that
in lemma \ref{section:moveC} (Fig. 4).
Next, because of the identity
$z\Gamma(z)=\Gamma(z+1)$, the quantity $\sup_{y\in{\R}}|\Gamma(x+iy)|$
decays superexponentially as $x\to-\infty$. Therefore,
by our estimate in corollary (\ref{section:estcorollary}) on the integral\\
$I_t(\eta,z)=\int_{X_R}\pi^*(\eta)(-\nabla_{tL}^2)^{-s}ds$,
(\ref{section:moveC2}) follows. \qed
\begin{figure}
\centerline {\psfig {figure=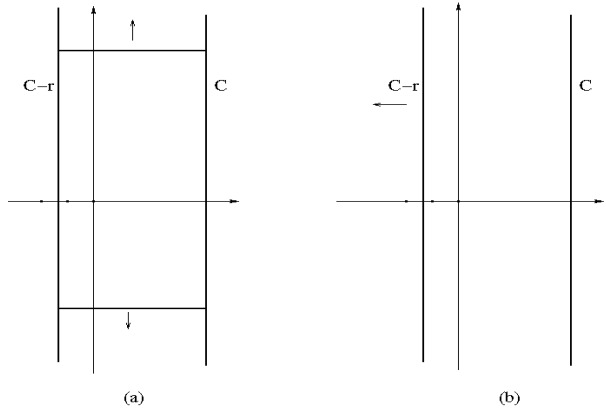,width=8cm, height=5.5cm}}
\caption{As $C$ is moved to the left, we pick up the residues.}
\end{figure}

Finally, we need to prove that as $R\to 0$,
only finitely many residues survive.
By theorem \ref{theorem:1},
$I_t(\eta,z)$ has only finitely many poles and they
are at most simple. Also, $\Gamma(z)$  has at most simple
poles, at $z=0,-1,-2\ldots$.
So, the residues of $\Gamma(z)I_t(\eta,z)$
are due to either one of these factors.
But the theorem \ref{theorem:1},
provides us with some knowledge of
the general algebraic form of $I_t(\eta,z)$.
It implies that
for $m>>0$, the residue at $z=-m$
will be a multiple of  a positive power of $R$,
so it will vanish as $R\to0$. \qed

\section{The de Rham operator on Riemanian surfaces.}

We consider the complexified vector bundle
$E\otimes{\C}=\Lambda^*M\otimes{\C}$ of exterior forms
over a compact manifold $M$ with no boundary.
The grading decomposition is that into differential forms of even
and odd degree. Given a Riemanian metric $g$ on $M,$ the
associated de Rham
operator $D=d+d^*$ has a well-known symbol $L=-\rho^2$.

In order to apply the theorem \ref{theorem:2}, we need to compute:

$$(-(\nabla+L)^2)^{-z}=(-(L^2+\nabla^2+[\nabla,L]))^{-z}.$$

The argument of the function $\nu\mapsto -\nu^{-z}$ can be
viewed as $-\rho^2$ plus some commuting perturbation which
is nilpotent. It follows that we may just use the Taylor
series expansion instead of the Cauchy integrals:

$$(-(\nabla+L)^2)^{-z}=\sum_{k=0}^{2n}(\rho^2)^{-(z+k)}
   \big(_{\,\,\,\,\,\,\,\,k}^{-(z+k)}\big)
    (\nabla^2+[\nabla,L])^{k},$$
where $\big(^z_k\big)$ means ${\frac {-z(-z-1)\ldots(-z-k+1)} {k!}}$.

For example, we consider the case when $M$ is a 2-manifold and
$\eta\equiv 1$.
Since that is exactly the todd class
for any 2-surface, by the Atiyah-Singer index theorem
both sides of theorem \ref{theorem:2} should give us
the euler characteristic. The left-hand side of theorem
\ref{theorem:2} yields:

\begin{align}\label{section:dR1}
\int_{T^*M}\exp((&\nabla+L)^2)=\\\nonumber
&=\int_{T^*M}\exp(-\rho^2)\exp(\nabla^2+[L,\nabla])\\\nonumber
&=\int_{T^*M}\exp(-\rho^2){\scriptstyle
    (1+(\nabla^2+[L,\nabla])+{\frac 1 2}(\nabla^2+[L,\nabla])^2+\ldots)}\\\nonumber
&=\int_{T^*M}\exp(-\rho^2){\scriptstyle ({\frac 1 {3!}}
         (\nabla^2+[L,\nabla])^3+{\frac 1 {4!}}(\nabla^2+[L,\nabla])^4)}
\end{align}
We keep only these two terms because they are the only ones that can
possibly contain a $4$-form which can be integrated over $T^*M.$
In fact, when we expand $(\nabla^2+[L,\nabla])^3$ we see that
only three terms really enter the picture, namely
$\nabla^2[L,\nabla]^2\,\,,$ $[L,\nabla]^2\nabla^2,$
and $[L,\nabla]\nabla^2[L,\nabla]$. From
$(\nabla^2+[L,\nabla])^4,$ the relevant term is $[L,\nabla]^4.$

{\lemma The term $[L,\nabla]^4$
vanishes as a section of $\Lambda^*T^*M\otimes End(\pi^*E)$, i.e.,
pointwise.\\}
{\bf Proof:} From section \ref{section:pt1},
$$[L,\nabla]=d_xL+d_\xi L+[\theta,L],$$
If one uses normal
coordinates on $M$ near some point $x$, then $\theta$,
 being comprised of Christoffel symbols is zero on
the fiber of $T^*M$ over $x$. The horizontal differential
$d_x L$ is also zero there.
Hence on that fiber, $[L,\nabla]=d_\xi L$ which is a matrix
of ``vertical'' forms on $T^*M$. Any power of it which is larger
than $dim(M)$ must vanish.\qed

Thus, the left-hand side of theorem \ref{theorem:2} reads:

\begin{equation}\label{section:chgen1}
\int_{T^*M}\exp(-\rho^2){\frac 1 6} {
(\nabla^2[L,\nabla]^2+[L,\nabla]^2\nabla^2+
                      [L,\nabla]\nabla^2[L,\nabla])}.
\end{equation}
Similar remarks apply on the right-hand side and we obtain:
\begin{align}\label{section:chgen2}
\lim_{R\to 0}\sum Res|_z\Gamma(z)\int_{X_R}
          {\frac {-z(z+1)(z+2)} 6}\rho^{-2(z+3)}\\\nonumber
   (\nabla^2[L,\nabla]^2+
                                     [L,\nabla]^2\nabla^2+
                                     [L,\nabla]\nabla^2[L,\nabla]).
\end{align}
Since ${\operatorname tr_s}{\scriptstyle   (\nabla^2[L,\nabla]^2+
                                     [L,\nabla]^2\nabla^2+
                                     [L,\nabla]\nabla^2[L,\nabla])}
                           =\omega\rho d\rho$
for some differential form $\omega$, it is enough to see that:
\begin{align}\label{section:chgen3}
\int_0^\infty \!\!\!\!exp(&-\rho^2)\rho d\rho\!=\\
\nonumber&=\lim_{R\to 0} \sum
Res|_z\Gamma(z)\int_R^\infty\!\!\!\!\!
           {(-z)(z+1)(z+2)} \rho^{-2(z+3)}\rho d\rho\nonumber\\
&={\frac 12}\lim_{R\to 0 }\sum Res|_z \Gamma(z+2)R^{-2(z+2)}
\nonumber\\\nonumber &={\frac 12}\lim_{R\to 0
}\sum_{m=0}^\infty{\frac {(-1)^m}{m!}}R^{2m}.
\end{align}

\section{The Chern character for a general spinor bundle.}

We apply our results to the Chern character of a spinor bundle
$S\to M$ associated to a vector bundle $\pi\co F\to M$, as
computed by Mathai and Quillen \cite{Q2}.  The role of
the cotangent bundle $\pi\co T^*M \to M$ is played by $F$
in this example, and the role of $E\to M$ is played by $S$.
So, theorem \ref{theorem:2}, strictly speaking does
not apply, though we could have proven it in a more general
setting. The reason we stated our theorem for $T^*M$
is that we have the Atiyah-Singer index theorem in mind,
for future applications.
Rather than applying theorem \ref{theorem:2}, we will
go through its proof. Namely,
we shall repeat  steps (2) and
(3) in a simpler way.

We proceed to outline the result of \cite{Q2}.
Some familiarity with spin structures is assumed
here. The reader can consult, e.g., {\it Spin Geometry}
by Lawson and Michelsohn for details \cite{LM}. We also warn that
the notation of \cite{Q2} is quite a bit different from
our own. We will briefly explain the differences in the end
of this section.

Let $\pi\co F\to M$ be a complex even-dimensional vector bundle with
a spin structure. In particular, this means that there is a
fiberwise metric $\mu$ on $F$. Let $S\to M$ be the associated
spinor bundle. The assumption of spin structure implies that $S$
can be split into a direct sum of even and odd subbundles:
$S=S^+\oplus S^-$, where the fibers $S_x^+$ and $S_x^-$ of $S^+$
and $S^-$ are the only two irreducible representations for the
spin group  of the fiber $F_x$. Thus, the spin structure induces a
${\Z}_2$-grading of $S$.

In order to form a Chern character we need a connection
$\nabla'$ on $S$ which respects that ${\Z}_2$-grading.
We also need an odd antiselfadjoint endomorphism $L$
of $\pi^*S$. The spinor bundle setup in \cite{Q2} requires,
among other things, that:
\begin{itemize}
\item $L$ be homogeneous of degree 1 in the radial
      coordinate $\rho$ of the fibers of $F$. Much as in the
      proof of theorems  \ref{theorem:1} and \ref{theorem:2},
      $\rho$ is induced by the metric $\mu$ and is given by
       $\rho(p)=\sqrt{\mu(p,p)}$ for all $p$ in $F$.
\item $\nabla'$ preserve the fiberwise metric $(\cdot,\cdot)$
      which is induced on $S$ by the metric $\mu$ of $F$.
      This means that for any two sections $\alpha$ and $\beta$
      of $S$, $$d(\alpha,\beta)=(\nabla'\alpha,\beta)+
                        (\alpha,\nabla'\beta).$$
\end{itemize}
The coordinate notation is the same as in  the proof of theorem
\ref{theorem:1}. The local coordinates on $F$ are the vertical
(fiberwise) cartesian coordinates are $\xi^1,\ldots, \xi^m$,
and the horizontal
coordinates $x_1,\ldots x_n$, which are also coordinates of $M$.
In fact, it makes sense to choose the $\mu$-orthonormal local frame
$e_1,\ldots e_m$ of $F$ and to choose coordinates
$\xi^j$ associated to that frame.
They may be replaced by spherical coordinates
$\rho$ and  $\Xi^1,\ldots,\Xi^{m-1}$ at our convenience.

To describe $L$ we recall that the spin structure of
$F$ stems from the fiberwise metric $\mu$. To begin with, we have the
{\it Clifford action} of $F$ on $S$ which is a fiberwise
${\R}$-linear bundle map $c\co F\to End(S),$ such that for any $(x,\xi)$
in $F_x$, $c(x,\xi)^2=-\mu_x(\xi,\xi)$. It is one of the standard
axioms for a Clifford actions that $c(x,\xi)$ be fiberwise
anti-selfadjoint endomorphism.
Thus,  each fiber $F_x$ is contained in a clifford algebra
$Cliff(F_x,\mu_x)$, which is a fiber of the bundle $Cliff(F,\mu)$.
Also, there is a map $$Cliff(F,\mu)\to End(S),$$ which is
an isomorphism of bundles and fiberwise an isomorphism
of algebras.
Now, the pullback $\pi^*F$ to the total space of $F$ is equipped with
the Clifford action on $\pi^*S$ which we shall also denote $c$
instead of $\pi^*c$. Let $\tau\co F\to\pi^*F$ be the tautological
section. Then the endomorphism $L=c(\tau(x,\xi))$ has all the required
properties. Its homogeneity in $\rho$ is obvious and
it is antiselfadjoint by hypothesis.

Abbreviating $c(\tau(e_j))$ by $\gamma_j$,  we may write
$L=\sum_j\xi^j\gamma_j$. Since the construction of the clifford
action on spinors using an orthonormal basis is completely
canonical, the coordinate expression for $L$ does not involve
the $x$-variables.  Since,
$${\frac12}(\gamma_i\gamma_j+\gamma_j\gamma_i)=-\delta_{ij},$$
it follows that $L^2=\mu(\xi,\xi)=-\rho^2$.

Next, to pick a connection on $\nabla'$, we start with a connection
on $F$ which is locally given by $d+\theta$. The connection $\nabla'$
on $S$ is constructed from it (see \cite{Q2} and \cite{LM}).
In order to describe the construction, we adopt the  {\it summation
notation}: we reserve the right to write any index as an
upper or a lower index. (Since we have chosen an orthonormal
local frame, there is no difference at all). Repetition of the
same index on the top and on the bottom implies summation.
Repetition on the top only or on the bottom only does {\bf not}.
The connection on $\pi^*S$ is given by
$$\nabla'=d+{\frac14}\theta^{ij}\gamma_i\gamma_j,$$
where $\theta^{ij}$ are just the matrix entries of the
endomorphism-valued$1$-form $\theta$. The connection
$\nabla=\pi^*\nabla'$ on $\pi^*S$
therefore makes sense. Observe that since $\gamma_j$ are odd,
the local endomorphism
$\theta^{ij}\gamma_i\gamma_j$ of $\pi^*S$ is even.
Moreover, it only depends on the variables $x_i$ and horizontal
differentials  $dx_i$, just as  before.
Therefore, the curvature of the connection $\nabla_L$
may be written as:
\begin {align}
\nabla_L^2&=\nabla^2+[\nabla,L]-\rho^2\\\nonumber
&=\nabla^2+d\xi^j\gamma_j+{\frac {\xi^k}4}
[\theta^{ij}\gamma_i\gamma_j,\gamma_k]-\rho^2\\\nonumber
&=\nabla^2+(d\xi^j)\gamma_j+{\frac {\xi^j}4}
\theta^{ij}\gamma_i-\rho^2
\end{align}
Here, just as in the previous example, the fact
that $L^2$ is a scalar is a tremendous simplification.
We may use Taylor series instead of Cauchy integrals and
we have the rule $e^{a+b}=e^ae^b.$

Now, the result from \cite{Q2} reads:
\begin{align}\label{section:QMch}
\operatorname{tr_s}
\exp\nabla^2_L&=
(-1)^{m/2}\Big({\frac i{2\pi }}\Big)^{-m/2}
\det\left({\frac{\sinh(\nabla^2/2)}{(\nabla^2/2)}}\right)^{\frac12}\\
&\,\,\qquad\times\operatorname{tr_s} ({\scriptstyle
\pi^{-m/2}e^{-\rho^2}\sum_{I}\varepsilon(I,I')
          Pf(\nabla^2/2)_I\prod_{j\in I'}
             ((d\xi_j)\gamma_j+{\frac {\xi^j}4}
                          \theta^i_j\gamma_i)}),
\nonumber
\end{align}
where:
\begin{itemize}
\item $I,I'$ are complementary (strictly increasing) multiindices
over the set $\{1,2,\ldots\,m\}$
and $\varepsilon(I,I')$ is a certain combinatorial $\pm 1$-valued
function of them, which shall not be relevant here.
\item $Pf(\nabla^2/2)_I$ is the Pfaffian of the submatrix
      of $\nabla^2/2$ determined by the multiindex $I$.
      For an unfamiliar reader, it suffices to know
      that it is a certain polynomial in the matrix
      entries of $\nabla^2$ which is just $1$ if
      $I$ is the empty multiindex.
\end{itemize}

For us, (\ref{section:QMch})  is greatly simplified by the fact that
we are only interested in the currents induced by
this Chern character on $\Omega^*M$.  Therefore, as
observed in theorem \ref{theorem:1}
and lemma \ref{theorem:poles}, we need  only
those terms of (\ref{section:QMch}) which involve {\it all}
the differentials $d\xi^j$, so the only relevant
multiindex is $I'=(1,2,\ldots,m)$, $I$ being
 empty and $Pf(\nabla^2/2)_I$ being $1$. The only
term of interest is therefore $$\pi^{-m/2}e^{-\rho^2}d\xi^1\ldots d\xi^m.$$
If we replace $L$ by $tL$, as in theorem \ref{theorem:2},
(\ref{section:QMch}) becomes
\begin{align}\label{section:QMchbis}
\operatorname{tr_s}
\exp\nabla^2_L&=
(-1)^{m/2}\Big({\frac i{2\pi }}\Big)^{-m/2}
\det\left({\frac{\sinh(\nabla^2/2)}{(\nabla^2/2)}}\right)^{\frac12}\\
&\qquad\qquad\qquad\times\operatorname{tr_s}
t^m\pi^{-m/2}e^{-t^2\rho^2}d\xi^1\ldots d\xi^m.\nonumber
\end{align}
Integrating this over any fiber of $F$, we  see that
for any $\eta\in\Omega^*M$,
\begin{align}\label{section:QMcurr}
\int_F\pi^*(\eta)\operatorname{tr_s}&
\exp\nabla^2_L=\\\nonumber&\int_F\pi^*(\eta) (-1)^{m/2}\Big({\frac
i{2\pi }}\Big)^{-m/2}
\det\left({\frac{\sinh(\nabla^2/2)}{(\nabla^2/2)}}\right)^{\frac12}.
\end{align}
This allows us to understand the residue formulation of this Chern
character.

A computation similar to the one in the previous example yields:
\begin{align}\label{section:st3QM4}
\exp\nabla^2_{L}&=
\sum_{k,l}{\frac1{2\pi i}}
          \int_{C}\!\!\!\!\big(^z_k\big)\Gamma(z)
          \rho^{-2(z+k)}
           {P_{k-l,l}(\nabla^2,[\nabla,L])}dz,
\end{align}
where by $P_{\mu,\nu}(A,B)$ we denote the homogeneous non-commutative
polynomial in $A$ and $B$ obtained by summing all the words
comprised of $\mu$ copies of $A$ and $\nu$ copies of $B$.

We now recall (\ref{section:QMch}). It involves the sum over multiindices
$I'$ and the only relevant multiindex was determined to be
$I'=(1,2\ldots\,,m)$, where $m$ is the fiberwise dimension of $F$.
This means that in (\ref{section:st3QM4}) only the terms with $l=m$
contribute to the current induced by Chern character on $\Omega^*M$.
We have seen a special case of this in the previous example,
where a normal coordinates argument was used to show that
only the terms which involve two copies of $[\nabla,L]$ are relevant.
(Recall from lemma \ref{theorem:poles} that such terms were called
{\it contributing}.)
In particular, it means that $[\nabla,L]$ contributes the vertical
differentials and no other differentials.

Coming back to our computation, the right-hand side of
(\ref{section:st3QM4}) is readily seen to be the Taylor
series for $\Gamma(z)(-\nabla^2_L)^{-2z}$. The discussion
in the previous paragraph implies that the contributing
part of $P_{k-l,l}(\nabla^2,[\nabla,L])$ is a multiple
of the vertical volume form:
$$d\xi^1\ldots d\xi^m=\rho^{m-1}d\rho d\Xi^1\ldots d\xi^{m-1}.$$
Thus, it supplies $m-1$ powers of $\rho$. It remains to
determine the residues, using the procedure from theorem
\ref{theorem:1}. Just as in that theorem, we set
$$X_R=_{_{def}}\{p\in F\co \rho(p)\ge R\},$$
and integrate from $R$ to $\infty$ with respect to $\rho$.
This, as we shall see, produces the residue at $(m-2k)/2$.
Let $\eta\in\Omega^\kappa M$, and express $(-\nabla^2_L)^{-z}$
as a sum of two differential forms:
$$(-\nabla^2_L)^{-z}=\nu_z+\omega_z d\rho,$$
where neither $\nu_z$ nor $\omega_z$ involve $d\rho$.
We obtain:
\begin{align}\label{section:st4QM1}
\Gamma(z)\int_{X_R}\pi^*(\eta)(-\nabla^2_L)^{-z} &=\int_R^\infty
\rho^{-2(z+k)+m-1}d\rho \int_{S^*M}\pi^*(\eta)\omega_z\\\nonumber
&={\frac {R^{-2(z+k)+m}}{2(z+k)+m}}
\int_{S^*M}\pi^*(\eta)\omega_z.
\end{align}
Counting the differential form degrees, we see that if
$\deg(\eta)=\kappa$, then the only contributing term
of (\ref{section:st3QM4}) is the $(k,m)$-th term. Here $k$
satisfies $2k-m=m+n-\kappa$. This term  produces a
residue at $m/2-k$ which is a current on $\kappa$-
forms. But according to (\ref{section:QMcurr}), the same
current is induced by the $(n-\kappa)$-component of
the differential form
$\det\left({\frac{\sinh(\nabla^2/2)}{(\nabla^2/2)}}\right)^{\frac12}$,
so that:
\begin{align}\int_M\eta\det&\left({\frac{\sinh(\nabla^2/2)}
                    {(\nabla^2/2)}}\right)^{\frac12}_{n-\kappa}=\nonumber\\
  &\lim_{R\to 0}Res|_{\frac{\kappa-m-n}2}
  \left(\Gamma(z)
  \int_{X_R}\pi^*(\eta)\operatorname{tr_s}(-\nabla_L^2)^{-z}\right),
\nonumber\end{align}
which agrees with lemma \ref{theorem:poles} if $m=n$.
\begin{remark}{\rm The condition $2k-m=m+n-\kappa$
  implies that we only have nonzero currents on $\kappa$
   forms if $\kappa $ is of the same parity as $n$.
   Thus, the location of the residue is an integer if $m$ is
   even and a half-integer if $m$ is odd.}
\end{remark}

{\bf Warning:}
   In \cite{Q2}, the relevant computation is in section 4, where
   $F$ is denoted by $E$, $\nabla^2$ is  denoted by $\Omega$ and
   the fiberwise coordinates $\xi^j$ are denoted by $x^j$.


\begin{thebibliography}{10}

\bibitem{Arfken}
G.~Arfken.
\newblock {\em Mathematical Methods for Physicists, 3rd ed.}
\newblock Academic Press, 1985.

\filbreak


\bibitem{AS0}
M.~F. Atiyah and I.~M. Singer.
\newblock The index of elliptic operators on compact manifolds.
\newblock {\em Bull. Amer. Math. Soc.}, 69:322--433, 1963.

\filbreak

\bibitem{AS123}
M.~F. Atiyah and I.~M. Singer.
\newblock The index of elliptic operators i, ii, iii.
\newblock {\em Ann. Math.}, 87:484--604, 1968.

\filbreak

\bibitem{AS4}
M.~F. Atiyah and I.~M. Singer.
\newblock The index of elliptic operators iv.
\newblock {\em Ann. Math.}, 93:119--138, 1971.

\filbreak

\bibitem{NCDG}
A.~Connes.
\newblock Noncommutative differential geometry.
\newblock {\em IHES Publ. Math.}, (62), 1985.

\filbreak

\bibitem{Con}
A.~Connes.
\newblock {\em Noncommutative geometry}.
\newblock Academic Press, 1994.

\filbreak

\bibitem{CM}
A.~Connes and H.~Moscovici.
\newblock The local index formula in noncommutative geometry.
\newblock {\em Geometric and Functional Analysis}, 5(2):174--243, 1995.

\filbreak


\bibitem{GVF}
J.~M. Gracia-Bondia; J. C. Varilly; H.~E. Figueroa.
\newblock {\em Elements of noncommutative geometry}.
\newblock Birkh\H auser Advanced Texts: Basler Lehrb\H ucher, 2001.

\filbreak


\filbreak


\bibitem{LM}
H.~B. Lawson and M.~Michelsohn.
\newblock {\em Spin Geometry}.
\newblock Princeton University Press, 1990.

\filbreak


\bibitem{Q2}
V.~Mathai and D.~Quillen.
\newblock Superconnections thom classes and equvariant differential forms.
\newblock {\em Topology}, 25(1):85--110, 1986.

\filbreak



\bibitem{Q3}
D.~Quillen.
\newblock Algebra cochains and cyclic cohomology.
\newblock {\em Publ. Math. I.H.E.S.}, 68:139--174, 1984.
\filbreak

\bibitem{Q}
D.~Quillen.
\newblock Superconnections and the chern character.
\newblock {\em Topology}, 24(1):89--95, 1985.

\filbreak




\end{thebibliography}
\end{document}